\documentclass[11pt,a4paper]{article}

\usepackage{graphicx}

\usepackage{amsmath}
\usepackage{verbatim}
\usepackage{amssymb}
\usepackage{hyperref}

\usepackage{amsthm}
\usepackage{geometry}
\usepackage{float}

\usepackage{setspace}
\spacing{1.5}
\usepackage{fancyhdr}

\usepackage{color} 

\usepackage{mathabx}
\usepackage{overpic}

\usepackage[all]{xy}

\graphicspath{{pics/}}

\newcommand{\al}{\alpha}

\newcommand{\vphi}{\varphi}

\newcommand{\be}{\beta}

\newcommand{\ga}{\gamma}

\newcommand{\de}{\delta}

\newcommand{\om}{\omega}

\newcommand{\na}{\nabla}

\newcommand{\NA}{\nabla}

\newcommand{\bs}{\boldsymbol}

\newcommand{\ra}{\rightarrow}

\newcommand{\lra}{\longrightarrow}

\newcommand{\Ra}{\Rightarrow}

\newcommand{\xra}{\xrightarrow}

\newcommand{\xlra}{\xlongrightarrow}

\newcommand{\rgl}{\rangle}

\newcommand{\lgl}{\langle}

\newcommand{\dash}{\textrm{-}}

\newcommand{\ot}{\otimes}

\newcommand{\bpf}{\begin{proof}}

\newcommand{\epf}{\end{proof}}

\newcommand{\bthm}{\begin{thm}}

\newcommand{\ethm}{\end{thm}}

\newcommand{\bprop}{\begin{prop}}

\newcommand{\eprop}{\end{prop}}

\newcommand{\bcor}{\begin{cor}}

\newcommand{\ecor}{\end{cor}}

\newcommand{\blem}{\begin{lem}}

\newcommand{\elem}{\end{lem}}

\newcommand{\bdefn}{\begin{defn}}

\newcommand{\edefn}{\end{defn}}

\newcommand{\bexmp}{\begin{exmp}}

\newcommand{\eexmp}{\end{exmp}}

\newcommand{\brem}{\begin{rem}}

\newcommand{\erem}{\end{rem}}

\newcommand{\bdia}{\begin{displaymath}\xymatrix}

\newcommand{\edia}{\end{displaymath}}

\newcommand{\beq}{\begin{equation*}\begin{aligned}}

\newcommand{\eeq}{\end{aligned}\end{equation*}}

\newcommand{\bref}{\textbf{Ref}}

\newcommand{\intg}{\mathbb{Z}}

\newcommand{\real}{\mathbb{R}}

\newcommand{\comp}{\mathbb{C}}

\newcommand{\quot}{\mathbb{H}}

\newcommand{\afv}{\mathbb{A}}

\newcommand{\prv}{\mathbb{P}}

\newcommand{\mco}{\mathcal{O}}

\newcommand{\mcc}{\mathcal{C}}

\newcommand{\mcf}{\mathcal{F}}

\newcommand{\mcg}{\mathcal{G}}

\newcommand{\mcs}{\mathcal{S}}

\newcommand{\cp}{\mathbb{CP}}

\newcommand{\mfo}{\mathfrak{O}}

\newcommand{\mfg}{\mathfrak{g}}

\newcommand{\msa}{\mathscr{A}}

\newcommand{\msr}{\mathscr{R}}

\newcommand{\msg}{\mathscr{G}}

\newcommand{\msd}{\mathscr{D}}

\newcommand{\itbf}{\item\textbf}

\newcommand{\seqa}{a_1,...,a_}

\newcommand{\seqx}{x_1,...,x_}

\newcommand{\seqy}{y_1,...,y_}

\newcommand{\seqf}{f_1,...,f_}

\newcommand{\cred}{\textcolor{red}}

\newcommand{\cblue}{\textcolor{blue}}

\newcommand{\mfa}{\mathfrak{a}}

\newcommand{\mfb}{\mathfrak{b}}

\newcommand{\mfm}{\mathfrak{m}}

\newcommand{\mfn}{\mathfrak{n}}

\newcommand{\mfp}{\mathfrak{p}}

\newcommand{\Af}{A_{(f)}}

\newcommand{\shm}{\underline{\rm SHM}}

\newtheorem{thm}{\textbf {Theorem}}[section]

\newtheorem{cor}[thm]{\textbf{Corollary}}

\newtheorem{prop}[thm]{\textbf{Proposition}}

\newtheorem{lem}[thm]{\textbf{Lemma}}

\newtheorem{conj}[thm]{Conjecture}

\theoremstyle{definition}

\newtheorem{defn}[thm]{\textbf{Definition}}

\newtheorem{exmp}[thm]{Example}

\newtheorem{notn}[thm]{Notation}

\theoremstyle{remark}

\newtheorem{rem}[thm]{Remark}

\def\cok{\operatorname{Coker}}

\newcommand{\txi}{\tilde{\xi}}

\newcommand{\bxi}{\bar{\xi}}
\usepackage{hyperref}
\usepackage{float}

\newcommand{\bz}{\bar{z}}

\DeclareMathOperator{\tr}{trunk}

\author{Nithin Kavi, Acton Boxborough Regional High School \\ Mentor: Zhenkun Li, MIT} 

\title{Cutting and Gluing Surfaces}
\date{}
\def\allfiles{}
\begin{document}
\bibliographystyle{plain}
\maketitle
\thispagestyle{empty}

\begin{abstract}
\thispagestyle{empty}
We start with a disk with $2n$ vertices along its boundary where pairs of vertices are connected with $n$ strips with certain restrictions. This forms a {\it pairing}. To relate two pairings, we define an operator called a cut-and-glue operation. We show that this operation does not change an invariant of pairings known as the {\it signature.} Pairings with a signature of $0$ are special because they are closely related to a topological construction through cut and glue operations that have other applications in topology. We prove that all balanced pairings for a fixed $n$ are connected on a surface with any number of boundary components. As a topological application, combined with works of Li, this shows that a properly embedded surface induces a well-defined grading on the sutured monopole Floer homology defined by Kronheimer and Mrowka. \end{abstract}

\tableofcontents

\section*{Acknowledgments}
I would like to thank my mentor Zhenkun Li for offering guidance throughout this project. Additionally, I thank the MIT math department and the MIT PRIMES program for giving me the opportunity to conduct this research.

\ifx\allfiles\undefined

\documentclass[12pt,a4paper]{article}


\usepackage{graphicx}

\usepackage{amsmath}

\usepackage{amssymb}
\DeclareMathOperator{\tr}{trunk}
\usepackage{amsthm}

\usepackage{geometry}

\usepackage{fancyhdr}

\usepackage{color} 









\newcommand{\al}{\alpha}

\newcommand{\vphi}{\varphi}

\newcommand{\be}{\beta}

\newcommand{\ga}{\gamma}

\newcommand{\de}{\delta}

\newcommand{\om}{\omega}

\newcommand{\na}{\nabla}

\newcommand{\NA}{\nabla}

\newcommand{\bs}{\boldsymbol}

\newcommand{\ra}{\rightarrow}

\newcommand{\lra}{\longrightarrow}

\newcommand{\Ra}{\Rightarrow}

\newcommand{\xra}{\xrightarrow}

\newcommand{\xlra}{\xlongrightarrow}

\newcommand{\rgl}{\rangle}

\newcommand{\lgl}{\langle}

\newcommand{\dash}{\textrm{-}}

\newcommand{\ot}{\otimes}

\newcommand{\bpf}{\begin{proof}}

\newcommand{\epf}{\end{proof}}

\newcommand{\bthm}{\begin{thm}}

\newcommand{\ethm}{\end{thm}}

\newcommand{\bprop}{\begin{prop}}

\newcommand{\eprop}{\end{prop}}

\newcommand{\bcor}{\begin{cor}}

\newcommand{\ecor}{\end{cor}}

\newcommand{\blem}{\begin{lem}}

\newcommand{\elem}{\end{lem}}

\newcommand{\bdefn}{\begin{defn}}

\newcommand{\edefn}{\end{defn}}

\newcommand{\bexmp}{\begin{exmp}}

\newcommand{\eexmp}{\end{exmp}}

\newcommand{\brem}{\begin{rem}}

\newcommand{\erem}{\end{rem}}

\newcommand{\bdia}{\begin{displaymath}\xymatrix}

\newcommand{\edia}{\end{displaymath}}

\newcommand{\beq}{\begin{equation*}\begin{aligned}}

\newcommand{\eeq}{\end{aligned}\end{equation*}}

\newcommand{\bref}{\textbf{Ref}}

\newcommand{\intg}{\mathbb{Z}}

\newcommand{\real}{\mathbb{R}}

\newcommand{\comp}{\mathbb{C}}

\newcommand{\quot}{\mathbb{H}}

\newcommand{\afv}{\mathbb{A}}

\newcommand{\prv}{\mathbb{P}}

\newcommand{\mco}{\mathcal{O}}

\newcommand{\mcc}{\mathcal{C}}

\newcommand{\mcf}{\mathcal{F}}

\newcommand{\mcg}{\mathcal{G}}

\newcommand{\mcs}{\mathcal{S}}

\newcommand{\cp}{\mathbb{CP}}

\newcommand{\mfo}{\mathfrak{O}}

\newcommand{\mfg}{\mathfrak{g}}

\newcommand{\msa}{\mathscr{A}}

\newcommand{\msr}{\mathscr{R}}

\newcommand{\msg}{\mathscr{G}}

\newcommand{\msd}{\mathscr{D}}

\newcommand{\itbf}{\item\textbf}

\newcommand{\seqa}{a_1,...,a_}

\newcommand{\seqx}{x_1,...,x_}

\newcommand{\seqy}{y_1,...,y_}

\newcommand{\seqf}{f_1,...,f_}

\newcommand{\cred}{\textcolor{red}}

\newcommand{\cblue}{\textcolor{blue}}

\newcommand{\mfa}{\mathfrak{a}}

\newcommand{\mfb}{\mathfrak{b}}

\newcommand{\mfm}{\mathfrak{m}}

\newcommand{\mfn}{\mathfrak{n}}

\newcommand{\mfp}{\mathfrak{p}}

\newcommand{\Af}{A_{(f)}}


\newtheorem{thm}{\textbf {Theorem}}[section]

\newtheorem{cor}[thm]{\textbf{Corollary}}

\newtheorem{prop}[thm]{\textbf{Proposition}}

\newtheorem{lem}[thm]{\textbf{Lemma}}

\newtheorem{conj}[thm]{Conjecture}

\newtheorem{conv}[thm]{Convention}

\newtheorem{prob}[thm]{Problem}

\newtheorem{exer}[thm]{Exercise}

\newtheorem{quest}[thm]{Question}

\theoremstyle{definition}

\newtheorem{defn}[thm]{\textbf{Definition}}

\newtheorem{defns}[thm]{Definitions}

\newtheorem{exmp}[thm]{Example}

\newtheorem{exmps}[thm]{Examples}

\newtheorem{var}[thm]{Variant}

\newtheorem{vars}[thm]{Variants}

\newtheorem{con}[thm]{Construction}

\newtheorem{notn}[thm]{Notation}

\newtheorem{notns}[thm]{Notations}

\theoremstyle{remark}

\newtheorem{rem}[thm]{Remark}

\newtheorem{rems}[thm]{Remarks}

\newtheorem{warn}[thm]{Warning}

\newtheorem{sch}[thm]{Scholium}

\newtheorem{expl}[thm]{Explanations}

\newtheorem*{theorem}{\textbf{Theorem}}

\newtheorem*{corollary}{\textbf{Corollary}}

\newtheorem*{proposition}{\textbf{Proposition}}

\newtheorem*{lemma}{\textbf{Lemma}}

\newtheorem*{example}{\textbf{Example}}

\def\cok{\operatorname{Coker}}

\newcommand{\txi}{\tilde{\xi}}

\newcommand{\bxi}{\bar{\xi}}

\newcommand{\bz}{\bar{z}}

\begin{document}

\bibliographystyle{plain}

\else

\fi

\section{Introduction}\label{sec_intro}


Balanced sutured manifolds and sutured monopole Floer homology are important tools in the study of $3$-dimensional topology. Balanced sutured manifolds were first introduced by Gabai \cite{gabai1983foliations} in 1983. They are a class of compact oriented $3$-manifolds $M$ together with a closed oriented $1$-submanifold $\ga\subset\partial M$, called the suture, which divides the boundary $\partial{M}$ into two parts of equal Euler characteristics. Later, in 2010, Kronheimer and Mrowka \cite{kronheimer2010knots} constructed the sutured monopole Floer homology on it. 

In 2019, Li \cite{li2019direct} constructed a grading on the monopole Floer homology associated to any properly embedded surface $S\subset M$ with connected boundary.

\bprop[Li, \cite{li2019direct}]\label{prop_grading_intro}
Suppose $(M,\ga)$ is a balanced sutured manifold and $S\subset M$ is an oriented properly embedded surface with connected boundary, then $S$ induces a well-defined grading on the sutured monopole Floer homology $SHM(M,\ga)$. 
\eprop

To construct the grading, Li first constructed a closed surface out of $S$, and the rest of the construction is straightforward. To construct the closed surface, he abstractly did the following: first he glued strips to $S$ along the intersection points of $S\cap\ga$, and then he assigned signs $\pm$ to each boundary components of the new surface in a fixed manner. If there was an equal number of positive boundary and negative boundary components, then he glued them in pairs and thus obtained a closed surface. 

This construction is the prototype of the balanced pairings we will discuss in the current paper. There are many different ways to glue the strips to $S$, and for a fixed $S\subset M$, one balanced pairing corresponds to a way to attach strips and thus a way to construct the grading on $SHM(M,\ga)$. So in order to prove the well-definedness of the grading associated to $S$, Li introduced the cut-and-glue operation to relate two balanced pairings.

\bprop[Li, \cite{li2019direct}]\label{prop_connected_pairings_intro}
If two balanced pairings are connected by a (finite) sequence of cut-and-glue operations, then they will result in the same grading on $SHM(M,\ga)$.
\eprop

In order to make the grading associated to $S$ well defined, Li only used a special class of balanced pairings to construct the closed surface, rather an arbitrary one; however, we prove one of his conjectures with the following theorem:

\begin{thm}\label{conbound}
For a surface with connected boundary, any two balanced pairings are connected by a sequence of cut-and-glue operations.
\end{thm}

Another important thing to notice is that Proposition \ref{prop_grading_intro} requires that $S$ has connected boundary. Clearly a general properly embedded surface $S\subset M$ can have more boundary components, so we also want to remove this extra condition to make the construction applicable to the general case. The definition of balanced pairings can be generalized naturally to the case where $S$ has multiple boundary components and Proposition \ref{prop_connected_pairings_intro} still holds. As another main result of the paper, we prove the following:

\bprop\label{multbound}
For any fixed surface, any two balanced pairings are connected by a sequence of cut-and-glue operations.
\eprop

\bcor\label{topapp}
Suppose $(M,\ga)$ is a balanced sutured manifold and $S\subset M$ is an oriented properly embedded surface. Then $S$ induces a well-defined grading on $SHM(M,\ga)$.
\ecor

In 2018, Li \cite{li2018gluing} used the grading to help compute the sutured monopole Floer homology of sutured solid torus as well as constructing minus versions of knot Floer homology. Hence we believe that Corollary \ref{topapp} could potentially be used to do more computations in monopole Floer homology and to construct minus versions for links in addition to knots (where Seifert surfaces of links might have more boundary components). 

The paper is organized as follows. In Section 2, we introduce basic definitions about our surfaces and show pictorially how we interpret the topological problem as a combinatorial one. In Section 3, we define a new invariant known as {\it signature}, which has not been defined previously, and prove its value does not change under a certain operation. 


In Section 4, we consider the case of balanced pairings where the signature is $0,$ and prove Theorem \ref{conbound}. We then extend this result in Section 5 to prove Proposition \ref{multbound}, which then implies Corollary \ref{topapp} by applying Proposition \ref{prop_connected_pairings_intro}. 

In Section 6, we summarize our results and provide possibilities for future research in this field of combinatorics.

Throughout this paper, we use color in figures to identify certain features of them more easily. However, the explanations that we provide can still be understood even if the figures are viewed in black and white.

\ifx\allfiles\undefined

\bibliography{Index}

\end{document}

\fi


\ifx\allfiles\undefined

\documentclass[12pt,a4paper]{article}


\usepackage{graphicx}

\usepackage{amsmath}

\usepackage{amssymb}
\DeclareMathOperator{\tr}{trunk}
\usepackage{amsthm}

\usepackage{geometry}

\usepackage{fancyhdr}

\usepackage{color} 









\newcommand{\al}{\alpha}

\newcommand{\vphi}{\varphi}

\newcommand{\be}{\beta}

\newcommand{\ga}{\gamma}

\newcommand{\de}{\delta}

\newcommand{\om}{\omega}

\newcommand{\na}{\nabla}

\newcommand{\NA}{\nabla}

\newcommand{\bs}{\boldsymbol}

\newcommand{\ra}{\rightarrow}

\newcommand{\lra}{\longrightarrow}

\newcommand{\Ra}{\Rightarrow}

\newcommand{\xra}{\xrightarrow}

\newcommand{\xlra}{\xlongrightarrow}

\newcommand{\rgl}{\rangle}

\newcommand{\lgl}{\langle}

\newcommand{\dash}{\textrm{-}}

\newcommand{\ot}{\otimes}

\newcommand{\bpf}{\begin{proof}}

\newcommand{\epf}{\end{proof}}

\newcommand{\bthm}{\begin{thm}}

\newcommand{\ethm}{\end{thm}}

\newcommand{\bprop}{\begin{prop}}

\newcommand{\eprop}{\end{prop}}

\newcommand{\bcor}{\begin{cor}}

\newcommand{\ecor}{\end{cor}}

\newcommand{\blem}{\begin{lem}}

\newcommand{\elem}{\end{lem}}

\newcommand{\bdefn}{\begin{defn}}

\newcommand{\edefn}{\end{defn}}

\newcommand{\bexmp}{\begin{exmp}}

\newcommand{\eexmp}{\end{exmp}}

\newcommand{\brem}{\begin{rem}}

\newcommand{\erem}{\end{rem}}

\newcommand{\bdia}{\begin{displaymath}\xymatrix}

\newcommand{\edia}{\end{displaymath}}

\newcommand{\beq}{\begin{equation*}\begin{aligned}}

\newcommand{\eeq}{\end{aligned}\end{equation*}}

\newcommand{\bref}{\textbf{Ref}}

\newcommand{\intg}{\mathbb{Z}}

\newcommand{\real}{\mathbb{R}}

\newcommand{\comp}{\mathbb{C}}

\newcommand{\quot}{\mathbb{H}}

\newcommand{\afv}{\mathbb{A}}

\newcommand{\prv}{\mathbb{P}}

\newcommand{\mco}{\mathcal{O}}

\newcommand{\mcc}{\mathcal{C}}

\newcommand{\mcf}{\mathcal{F}}

\newcommand{\mcg}{\mathcal{G}}

\newcommand{\mcs}{\mathcal{S}}

\newcommand{\cp}{\mathbb{CP}}

\newcommand{\mfo}{\mathfrak{O}}

\newcommand{\mfg}{\mathfrak{g}}

\newcommand{\msa}{\mathscr{A}}

\newcommand{\msr}{\mathscr{R}}

\newcommand{\msg}{\mathscr{G}}

\newcommand{\msd}{\mathscr{D}}

\newcommand{\itbf}{\item\textbf}

\newcommand{\seqa}{a_1,...,a_}

\newcommand{\seqx}{x_1,...,x_}

\newcommand{\seqy}{y_1,...,y_}

\newcommand{\seqf}{f_1,...,f_}

\newcommand{\cred}{\textcolor{red}}

\newcommand{\cblue}{\textcolor{blue}}

\newcommand{\mfa}{\mathfrak{a}}

\newcommand{\mfb}{\mathfrak{b}}

\newcommand{\mfm}{\mathfrak{m}}

\newcommand{\mfn}{\mathfrak{n}}

\newcommand{\mfp}{\mathfrak{p}}

\newcommand{\Af}{A_{(f)}}


\newtheorem{thm}{\textbf {Theorem}}[section]

\newtheorem{cor}[thm]{\textbf{Corollary}}

\newtheorem{prop}[thm]{\textbf{Proposition}}

\newtheorem{lem}[thm]{\textbf{Lemma}}

\newtheorem{conj}[thm]{Conjecture}

\newtheorem{conv}[thm]{Convention}

\newtheorem{prob}[thm]{Problem}

\newtheorem{exer}[thm]{Exercise}

\newtheorem{quest}[thm]{Question}

\theoremstyle{definition}

\newtheorem{defn}[thm]{\textbf{Definition}}

\newtheorem{defns}[thm]{Definitions}

\newtheorem{exmp}[thm]{Example}

\newtheorem{exmps}[thm]{Examples}

\newtheorem{var}[thm]{Variant}

\newtheorem{vars}[thm]{Variants}

\newtheorem{con}[thm]{Construction}

\newtheorem{notn}[thm]{Notation}

\newtheorem{notns}[thm]{Notations}

\theoremstyle{remark}

\newtheorem{rem}[thm]{Remark}

\newtheorem{rems}[thm]{Remarks}

\newtheorem{warn}[thm]{Warning}

\newtheorem{sch}[thm]{Scholium}

\newtheorem{expl}[thm]{Explanations}

\newtheorem*{theorem}{\textbf{Theorem}}

\newtheorem*{corollary}{\textbf{Corollary}}

\newtheorem*{proposition}{\textbf{Proposition}}

\newtheorem*{lemma}{\textbf{Lemma}}

\newtheorem*{example}{\textbf{Example}}

\def\cok{\operatorname{Coker}}

\newcommand{\txi}{\tilde{\xi}}

\newcommand{\bxi}{\bar{\xi}}

\newcommand{\bz}{\bar{z}}



\begin{document}

\bibliographystyle{plain}

\else

\fi

\section{Preliminaries}
We will start with some necessary definitions.

\begin{figure}[H]
    \centering
    \includegraphics[scale = .175]{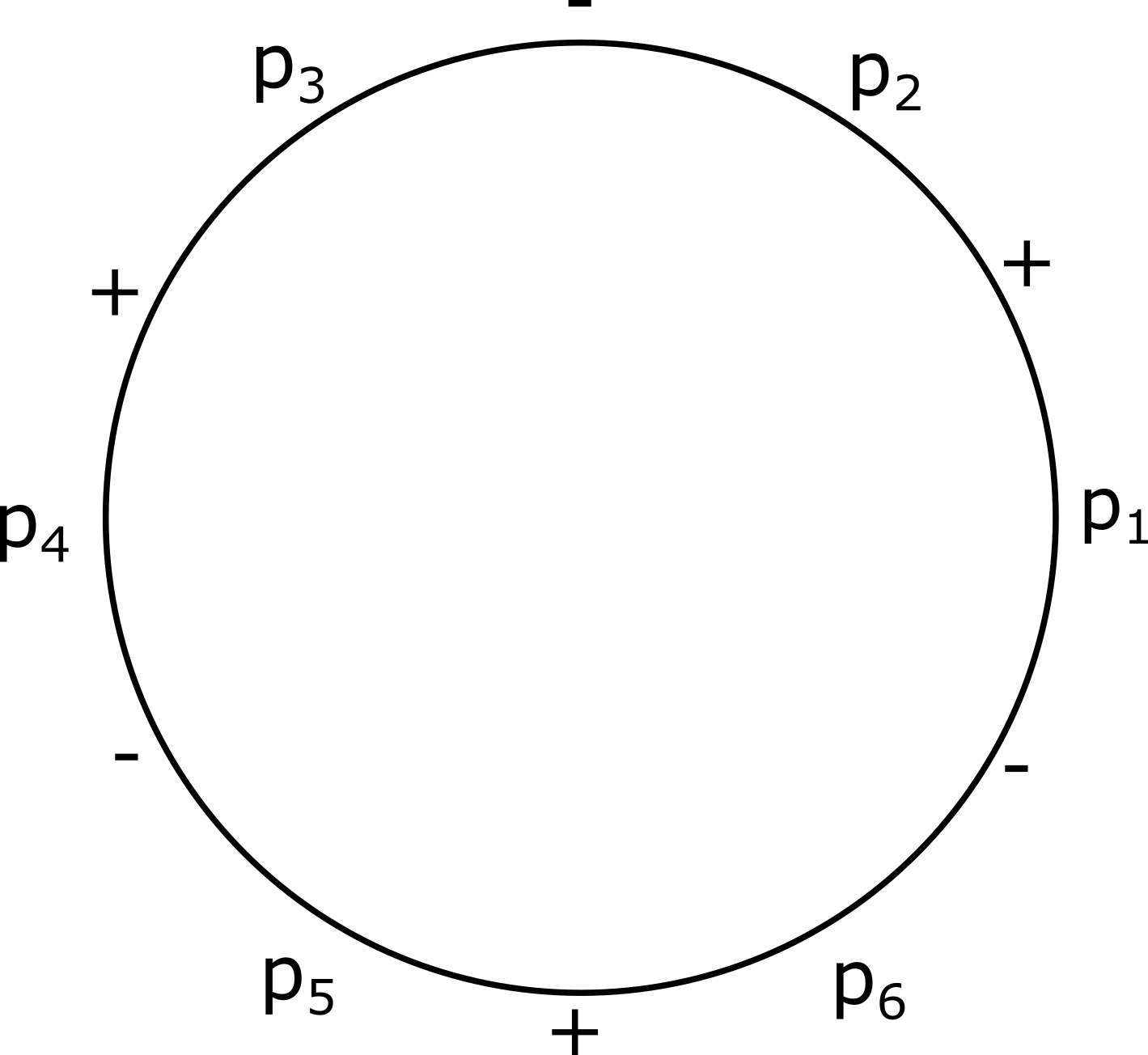}
    \caption{An example of the vertices on the disk.}
    \label{fig1}
\end{figure}

We begin by considering the boundary $S^1$ of a disk $D$. Regard $S^1$ as the unit circle on the complex plane $\mathbb{C}$, and let the points $p_1, p_2, \ldots, p_{2n}$ be on $S^1$ such that $p_k$ is the point $e^{2k \pi i/n}.$ In addition, we require $n$ to be odd. We have $2n$ arcs on the boundary $S^1$ which are separated by vertices $p_1,...,p_{2n}$. We assign $+$ and $-$ to each arc alternately, where we arbitrarily require that the arc between $p_1$ and $p_2$ is positive. Then, it follows that the arc between $p_2$ and $p_3$ is negative, the arc between $p_3$ and $p_4$ is positive, and so on. We call these arcs either plus-arcs or minus-arcs according to the sign. An example of this where $n = 3$ is shown by Figure \ref{fig1}.

We attach $n$ strips $S_1,...,S_n$ onto the disk $D$ along the vertices. Each strip is topologically an $[-\varepsilon,\varepsilon]\times[-1,1]$. For each $i=1,...,n$, we pick two vertices $p_{j_i}$, $p_{k_i}$ so that $j_i$ and $k_i$ are of different parity, and glue the strip $S_i$ to $D$ via an embedding
$$\rho_i:[-\varepsilon,\varepsilon]\times\{\pm1\}\ra S^1=\partial{D}.$$
The embedding $\rho_i$ should satisfy the following restrictions:

(1). We shall require that $\rho(\{0\}\times\{-1\})=p_{j_i}$.

(2). We shall require that $\rho(\{0\}\times\{1\})=p_{k_i}$.

(3). We shall require that $\rho(\{-\varepsilon\}\times\{\pm1\})$ lie in minus-arcs on $S^1$ while $\rho(\{+\varepsilon\}\times\{\pm1\})$ lie in plus-arcs on $S^1$.

We shall require that the $n$ strips are all attached along disjoint vertices.

\begin{notn}\label{not_1}
For later convenience, if a strip is identified with $S_i=[-\varepsilon,\varepsilon]\times[-1,1]$, we write $\partial_-S_i=\{-\varepsilon\}\times[-1,1]$ and write $\partial_+S_i=\{-\varepsilon\}\times[-1,1]$.
\end{notn}

Let $\tilde{D}=D\cup S_1\cup...\cup S_n$. Then, recall that $\chi(D) = 1$ (here $\chi$ refers to the Euler characteristic). It is well known that adding a strip to a surface decreases its Euler characteristic by $1$, so after adding the $n$ strips $S_1, S_2, \ldots, S_n$, we get $\chi(\tilde{D}) = 1 - n.$ 

\begin{lem}\label{evenbd} The number of boundary components of $\tilde{D}$ is even.

\bpf

Suppose $\tilde{D}$ has $d$ boundary components. Then we get $g(\tilde{D}) = \frac{n + 1 - d}{2},$ where $g(\tilde{D})$ is the genus of $\tilde{D}$. Since $n$ is odd and the genus is always an integer, we get that $d$ is even. 

\epf 

\end{lem}

\begin{figure}[H]
    \centering
    \includegraphics[scale=0.2]{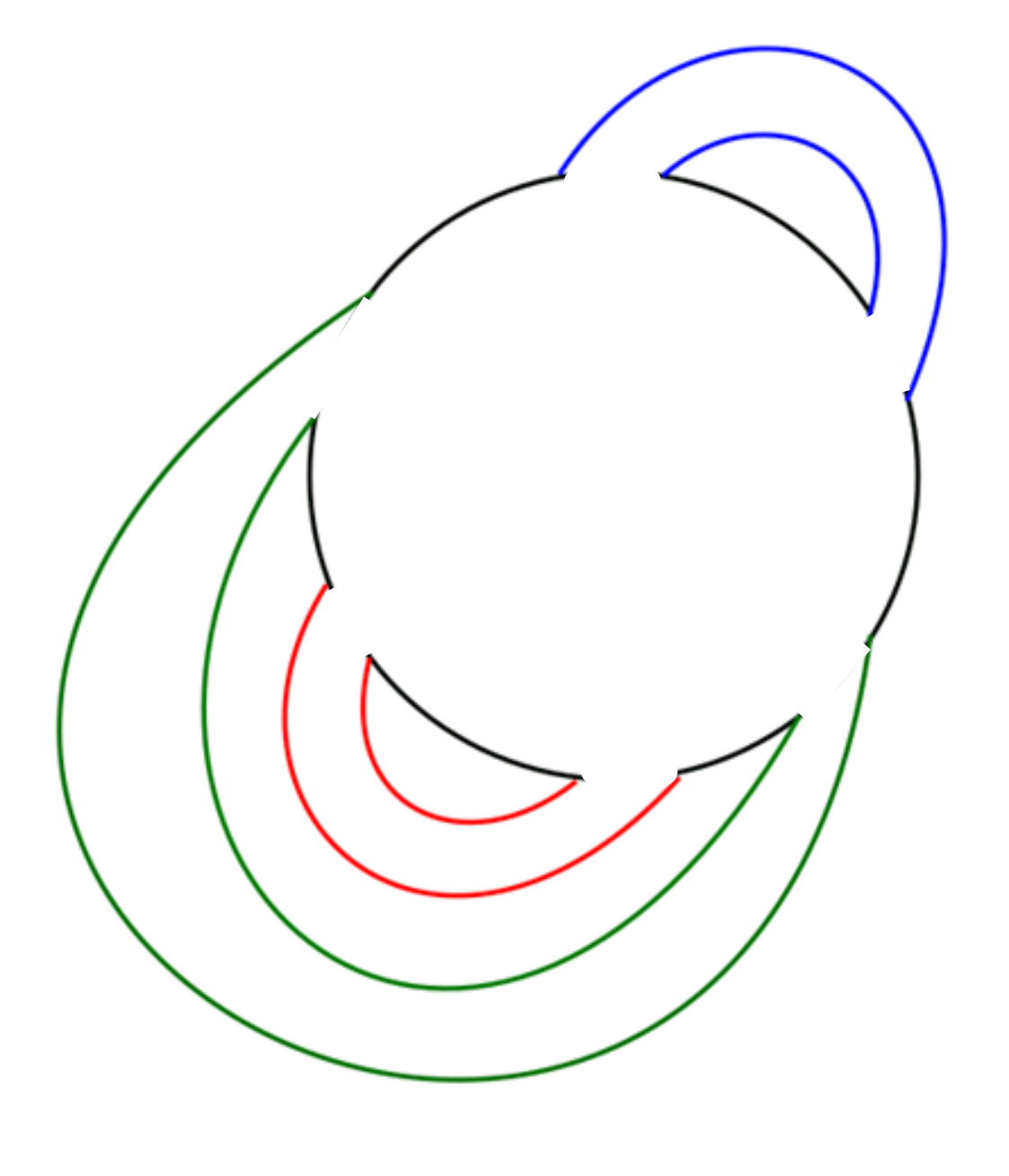}
    \caption{Example of $D_\mathcal{P}$ embedded in $\mathbb{C}$ for $n = 3$.}
    \label{fig:my_intersection}
\end{figure}

\begin{defn}\label{pairingdefn}
A \textit{pairing} is a set of $n$ couples $\{(i_1,j_1),...,(i_n,j_n)\}$, so that

(1). For $k=1,...,n$, $i_k\equiv j_k+1~({\rm mod}~2)$.

(2). We have $\{i_1,j_1,...,i_n,j_n\}=\{1,2,...,2n\}$.

Here $n$ is an odd positive integer and is called the {\it size} of the pairing. We use $\mathcal{P}_n$ to denote a pairing of size $n$, and use $\Pi_n$ to denote the set of all pairings of size $n$. When $n$ is clear, we might write just $\mathcal{P}$ and $\Pi$.

A pairing $\mathcal{P}$ gives us a unique way to attach $n$ strips to the disk $D$ and we call the resulting surface $\tilde{D}_{\mathcal{P}}$. See Figure \ref{fig:my_intersection} for an example of $\tilde{D}_{\mathcal{P}}$ when $n = 3.$ 
\end{defn}

\begin{defn}\label{posneg} A boundary component $\al$ of $\tilde{D}_{\mathcal{P}}$ is called \textit{positive} if $\al\cap {D}$ consists of only plus-arcs. It is called \textit{negative} if $\al\cap {D}$ consists of only minus-arcs.
\end{defn}

\begin{rem} Note that we never have a boundary component of $\tilde{D}_{\mathcal{P}}$ with both plus-arcs and minus-arcs because of the way we attach the strips to the disk $D.$ This is why requirement (1) in Definition \ref{pairingdefn} exists.

\end{rem}


\ifx\allfiles\undefined

\bibliography{Index}

\end{document}

\fi

\ifx\allfiles\undefined

\documentclass[12pt,a4paper]{article}


\usepackage{graphicx}

\usepackage{amsmath}

\usepackage{amssymb}
\DeclareMathOperator{\tr}{trunk}
\usepackage{amsthm}

\usepackage{geometry}

\usepackage{fancyhdr}

\usepackage{color} 

\usepackage{mathabx}

\usepackage{overpic}









\newcommand{\shm}{\underline{\rm SHM}}

\newcommand{\al}{\alpha}

\newcommand{\vphi}{\varphi}

\newcommand{\be}{\beta}

\newcommand{\ga}{\gamma}

\newcommand{\de}{\delta}

\newcommand{\om}{\omega}

\newcommand{\na}{\nabla}

\newcommand{\NA}{\nabla}

\newcommand{\bs}{\boldsymbol}

\newcommand{\ra}{\rightarrow}

\newcommand{\lra}{\longrightarrow}

\newcommand{\Ra}{\Rightarrow}

\newcommand{\xra}{\xrightarrow}

\newcommand{\xlra}{\xlongrightarrow}

\newcommand{\rgl}{\rangle}

\newcommand{\lgl}{\langle}

\newcommand{\dash}{\textrm{-}}

\newcommand{\ot}{\otimes}

\newcommand{\bpf}{\begin{proof}}

\newcommand{\epf}{\end{proof}}

\newcommand{\bthm}{\begin{thm}}

\newcommand{\ethm}{\end{thm}}

\newcommand{\bprop}{\begin{prop}}

\newcommand{\eprop}{\end{prop}}

\newcommand{\bcor}{\begin{cor}}

\newcommand{\ecor}{\end{cor}}

\newcommand{\blem}{\begin{lem}}

\newcommand{\elem}{\end{lem}}

\newcommand{\bdefn}{\begin{defn}}

\newcommand{\edefn}{\end{defn}}

\newcommand{\bexmp}{\begin{exmp}}

\newcommand{\eexmp}{\end{exmp}}

\newcommand{\brem}{\begin{rem}}

\newcommand{\erem}{\end{rem}}

\newcommand{\bdia}{\begin{displaymath}\xymatrix}

\newcommand{\edia}{\end{displaymath}}

\newcommand{\beq}{\begin{equation*}\begin{aligned}}

\newcommand{\eeq}{\end{aligned}\end{equation*}}

\newcommand{\bref}{\textbf{Ref}}

\newcommand{\intg}{\mathbb{Z}}

\newcommand{\real}{\mathbb{R}}

\newcommand{\comp}{\mathbb{C}}

\newcommand{\quot}{\mathbb{H}}

\newcommand{\afv}{\mathbb{A}}

\newcommand{\prv}{\mathbb{P}}

\newcommand{\mco}{\mathcal{O}}

\newcommand{\mcc}{\mathcal{C}}

\newcommand{\mcf}{\mathcal{F}}

\newcommand{\mcg}{\mathcal{G}}

\newcommand{\mcs}{\mathcal{S}}

\newcommand{\cp}{\mathbb{CP}}

\newcommand{\mfo}{\mathfrak{O}}

\newcommand{\mfg}{\mathfrak{g}}

\newcommand{\msa}{\mathscr{A}}

\newcommand{\msr}{\mathscr{R}}

\newcommand{\msg}{\mathscr{G}}

\newcommand{\msd}{\mathscr{D}}

\newcommand{\itbf}{\item\textbf}

\newcommand{\seqa}{a_1,...,a_}

\newcommand{\seqx}{x_1,...,x_}

\newcommand{\seqy}{y_1,...,y_}

\newcommand{\seqf}{f_1,...,f_}

\newcommand{\cred}{\textcolor{red}}

\newcommand{\cblue}{\textcolor{blue}}

\newcommand{\mfa}{\mathfrak{a}}

\newcommand{\mfb}{\mathfrak{b}}

\newcommand{\mfm}{\mathfrak{m}}

\newcommand{\mfn}{\mathfrak{n}}

\newcommand{\mfp}{\mathfrak{p}}

\newcommand{\Af}{A_{(f)}}


\newtheorem{thm}{\textbf {Theorem}}[section]

\newtheorem{cor}[thm]{\textbf{Corollary}}

\newtheorem{prop}[thm]{\textbf{Proposition}}

\newtheorem{lem}[thm]{\textbf{Lemma}}

\newtheorem{conj}[thm]{Conjecture}

\newtheorem{conv}[thm]{Convention}

\newtheorem{prob}[thm]{Problem}

\newtheorem{exer}[thm]{Exercise}

\newtheorem{quest}[thm]{Question}

\theoremstyle{definition}

\newtheorem{defn}[thm]{\textbf{Definition}}

\newtheorem{defns}[thm]{Definitions}

\newtheorem{exmp}[thm]{Example}

\newtheorem{exmps}[thm]{Examples}

\newtheorem{var}[thm]{Variant}

\newtheorem{vars}[thm]{Variants}

\newtheorem{con}[thm]{Construction}

\newtheorem{notn}[thm]{Notation}

\newtheorem{notns}[thm]{Notations}

\theoremstyle{remark}

\newtheorem{rem}[thm]{Remark}

\newtheorem{rems}[thm]{Remarks}

\newtheorem{warn}[thm]{Warning}

\newtheorem{sch}[thm]{Scholium}

\newtheorem{expl}[thm]{Explanations}

\newtheorem*{theorem}{\textbf{Theorem}}

\newtheorem*{corollary}{\textbf{Corollary}}

\newtheorem*{proposition}{\textbf{Proposition}}

\newtheorem*{lemma}{\textbf{Lemma}}

\newtheorem*{example}{\textbf{Example}}

\def\cok{\operatorname{Coker}}

\newcommand{\txi}{\tilde{\xi}}

\newcommand{\bxi}{\bar{\xi}}

\newcommand{\bz}{\bar{z}}

\begin{document}

\bibliographystyle{plain}

\else

\fi
\section{Signature}

Recall that we have a disk $D$ with $2n$ vertices. If $\mathcal{P}$ is a pairing then we can attach $n$ strips to $D$ to get a surface $\tilde{D}_{\mathcal{P}}$. 

\begin{defn} The \textit{signature} of a pairing ${\mathcal{P}}$, which we denote by $\sigma({\mathcal{P}})$, is the difference of the number of positive boundary components on $\tilde{D}_{\mathcal{P}}$ and the number of negative boundary components on $\tilde{D}_{\mathcal{P}}$.
\end{defn}

\brem Since the total number of boundary components of $\tilde{D}$ is even by Lemma \ref{evenbd}, the number of positive and negative boundary components must have the same parity. This means that the signature is always an even integer. 

\erem 

\begin{defn}\label{balanced} A pairing ${\mathcal{P}}$ is \textit{balanced} if $\sigma({\mathcal{P}})=0.$

\end{defn}

\begin{lem}\label{special} Let $\mathcal{P}_s$ be the pairing
$$\mathcal{P}_s=\{(1,n+1),(2,n+2),...,(n,2n)\},$$
Then $P_s$ is balanced for any $n.$ In particular, $P_s$ always has exactly one positive boundary and exactly one negative boundary.
\bpf
Note that since $n$ is odd, $k$ and $n + k$ always have opposite parity as required. Begin at the number $1$ and go along the plus-arc to vertex $2$. Then $2$ is connected to $n + 2,$ and then you go to $n + 3,$ and then to $3$, and then to $4$ and then to $n + 4$ etc, going around the circle until you reach the number $n + 1.$ At this point you go back to $1$ and this is the only positive boundary as it includes all the vertices on the circle. The negative boundary is analogous. We also note that due to the rotational symmetry of $P_s,$ it is independent of which vertex we label as $p_1.$
\epf
\end{lem}

\begin{exmp} Figure \ref{fig:special} shows $P_s$ for $n = 3.$

\begin{figure}[H]
    \centering
    \includegraphics[scale = 0.25]{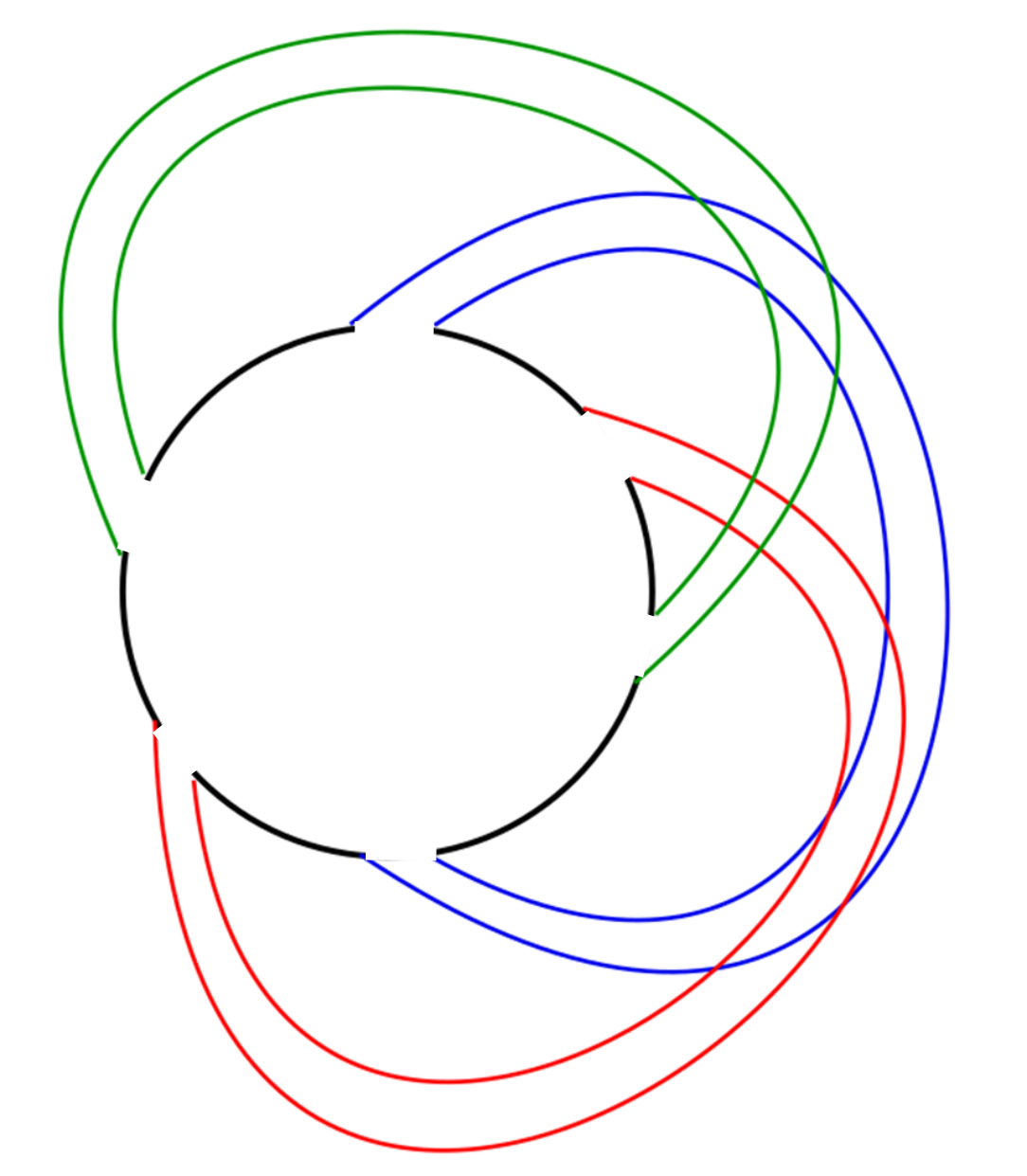}
    \caption{$D_{\mathcal{P}}$ has one positive boundary and one negative boundary.}
    \label{fig:special}
\end{figure}
\end{exmp}

\begin{defn} Suppose $n$ is a fixed odd integer and $i, j$ are integers so that $1 \leq i, j \leq 2n$. Then we define a function
$$\tau_i:\{1,...,2n\}\ra\{1,...,2n\}$$
as follows: if $i \leq j,$ then we let $\tau_i(j) = j - i + 1,$ and if $i > j,$ then we let $\tau_i(j) = 2n - (i - j) + 1.$

\end{defn}

\begin{rem}

The $\tau$ function allows us to effectively renumber the vertices on the disk beginning with relabeling $i$ as $1.$ This allows us to understand the order of the numbers with respect to a specific number $i.$

\end{rem} 

\begin{defn} For a fixed pairing $\mathcal{P}$, suppose we have $p_i, p_j, p_k, p_l$ on the surface of the disk where $p_i$ is connected to $p_k$ by a strip and $p_j$ is connected to $p_l$ by another strip. If we have 
\begin{equation}\label{eq_1}
    \tau_i(j) < \tau_i(k) < \tau_i(l),
\end{equation}
then $\mathcal{P}$ is called \textit{non-planar}. If for any $p_i, p_j, p_k, p_l$ where $p_i$ is connected to $p_k$ and $p_j$ is connected to $p_l$ the inequality (\ref{eq_1}) is not true, then $\mathcal{P}$ is a \textit{planar} pairing.
\end{defn}

\brem
A more geometric explanation of whether a pairing is planar:
\erem

\begin{figure}[H]
    \centering
    \includegraphics[scale = 0.25]{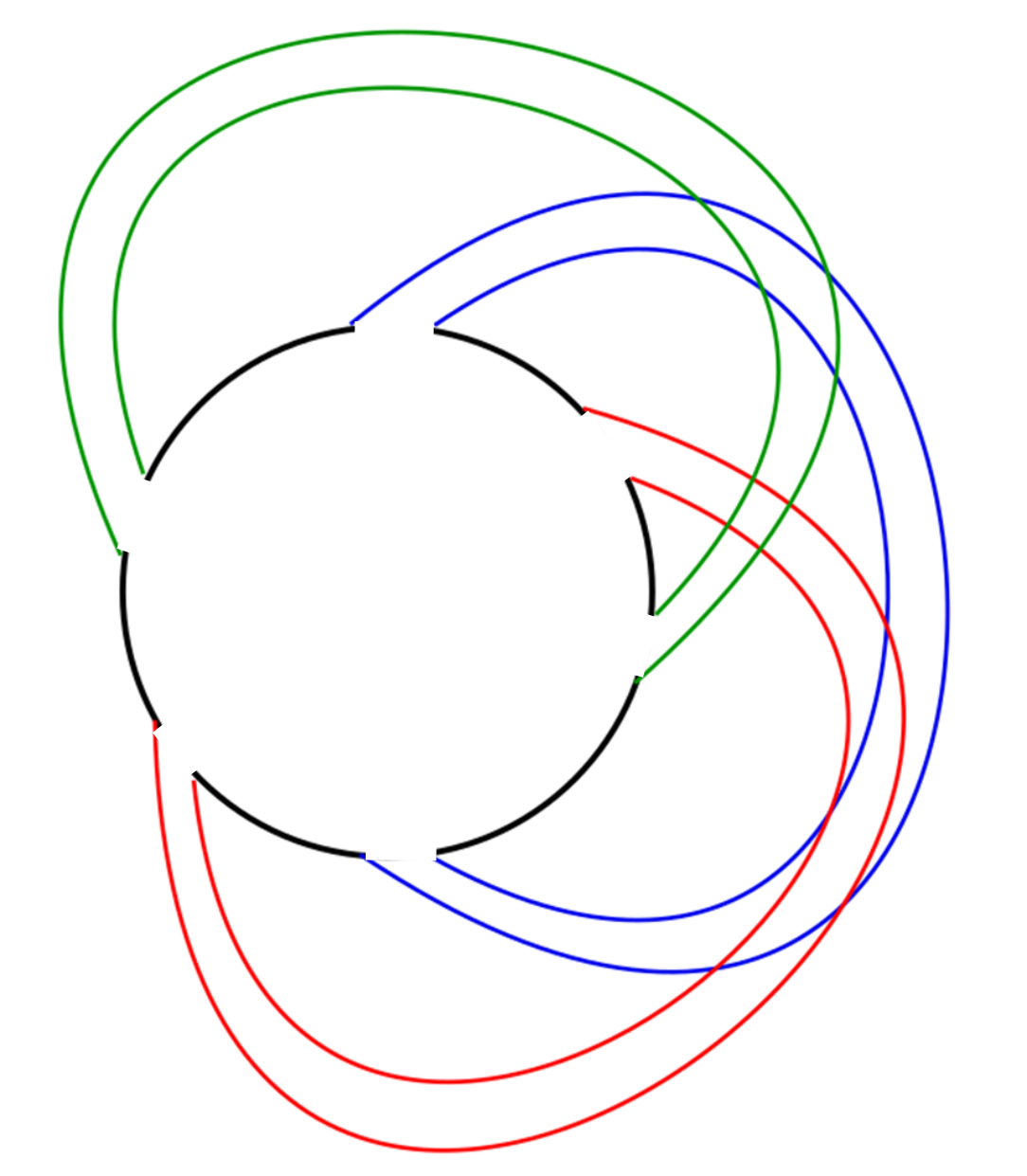}
    \caption{$\mathcal{P}_s$ as defined in Definition \ref{special} is nonplanar.}
    \label{fig4}
\end{figure}

The strips ``overlap" in the geometrical 2-dimensional embedding shown in Figure \ref{fig4}, so this pairing is nonplanar.

\begin{defn}\label{cutglue} Suppose $\mathcal{P}$ is a pairing and that with respect to this pairing $\mathcal{P}$, vertices $p_{e_1}$ and $p_{o_1}$ are connected by a strip $S_1$, and vertices $p_{e_2}$ and $p_{o_2}$ are connected by a strip $S_2$. Here the $e_i$ are even and the $o_i$ are odd. If the four arcs $\partial{\pm}S_1$ and $\partial{\pm}S_2$ are on either exactly $2$ or $4$ distinct boundary components of $\tilde{D}_{\mathcal{P}}$, then we can form a new pairing $\mathcal{P}'$ by connecting $e_1$ and $o_2$ together with a strip and connecting $e_2$ and $o_1$ together with a strip while leaving the other pairs unchanged. We call this process a {\it cut-and-glue operation}.

If the four arcs are on exactly $3$ distinct boundary components of $\tilde{D}_{\mathcal{P}}$, then no such cut-and-glue operation is allowed.
\end{defn}

A more geometric interpretation of the cut-and-glue operation is as follows. Suppose for $i=1,2$, the strip $S_i$ is parameterized by $S_i=[-\varepsilon,\varepsilon]\times [-1,1]_i$, and is glued to the disk $D$ via 
$$\rho_i:[-\varepsilon,\varepsilon]\times\{\pm1\}\ra\partial D.$$
We shall require that $\rho_i(\{0\}\times\{-1\})=p_{e_i}$ and $\rho_i(\{0\}\times\{1\})=p_{o_i}$. The cut-and-glue operation corresponds to replacing $S_1$ and $S_2$ by the new strips $S_1'$ and $S_2'$ as follows: we can cut $S_i=[-\varepsilon,\varepsilon]\times [-1,1]_i$ open along the arc $[-\varepsilon,\varepsilon]\times\{0\}_i$, and let $S_{i,-}=[-\varepsilon,\varepsilon]\times [-1,0]_i$ and $S_{i,+}=[-\varepsilon,\varepsilon]\times [0,1]_i$. Then we reglue $S_{1,-}$ to $S_{2,+}$ along the identity on $[-\varepsilon,\varepsilon]\times\{0\}$, and reglue $S_{1,+}$ to $S_{2,-}$ along the identity on $[-\varepsilon,\varepsilon]\times\{0\}$.

\begin{prop}\label{consig} A cut-and-glue operation on a pairing does not change its signature.

\end{prop}

\begin{proof} 
Suppose the cut-and-glue operation involves four vertices and two strips $S_1, S_2$ as in Definition \ref{cutglue}. Then we have two cases.

Case 1: The four arcs are in two separate boundary components $\be_+$ and $\be_-$. Note $\be_+$ contains only plus-arcs and $\be_-$ contains only minus-arcs as in Definition \ref{posneg}.

We can represent $\be_+$ by $b_1b_2b_3 \ldots b_k,$ where the $b_i$ denote either a plus-arc or $\partial_+{S_i}$ for some strip $S_i$. Similarly, $\be_-$ can be represented by $c_1c_2c_3 \ldots c_l$ where the $c_i$ are analogous to the $b_i.$ Suppose $b_p=\partial_+S_1$ and $b_q=\partial_+{S}_2$, similarly assume $c_r=\partial_-S_1$ and $c_s=\partial_-{S}_2$. Let us also assume that $p<q$ and $r<s$. Recall that in the cut-and-glue operation, we cut $S_1$ and $S_2$ in the middle and thus $b_p$, $b_q$, $c_r$ and $c_s$ are cut into eight pieces $b_{p,\pm}$, $b_{q,\pm}$ $c_{r,\pm}$ and $c_{s,\pm}$.

After the re-gluing, the arcs are re-grouped as $b_{p,+}b_{q,-}$ $b_{q,+}b_{p,-}$, $c_{r,+}c_{s,-}$ and $c_{s,+}c_{r,-}$. Then we get two positive boundaries: 

$$b_1b_2 \ldots b_{p-1}(b_{p,+}b_{q,-})b_{q+1} \ldots b_{k}, \hspace{5mm}b_{p+2} \ldots b_{q-1}(b_{q,+}b_{p,-}b_{p+1}).$$ Similarly, we get two negative boundaries: 
$$c_1c_2 \ldots c_{s - 1}(c_{s,-}c_{r+})c_{r + 1} \ldots c_{l}, \hspace{5mm} c_{s + 2} \ldots c_{r,-}c_{s,+}c_{s + 1}.$$ Thus the number of positive boundaries and the number of negative boundaries both increase by $1,$ so the signature remains the same.

Case 2: The four arcs are in four separate boundary components.

This case is the reverse process of Case 1, so the number of positive boundaries and the number of negative boundaries both decrease by $1$. The signature remains the same in both cases as desired.

\end{proof}








%

\ifx\allfiles\undefined
\end{document}

\fi

\ifx\allfiles\undefined

\documentclass[12pt,a4paper]{article}


\usepackage{graphicx}

\usepackage{amsmath}

\usepackage{amssymb}
\DeclareMathOperator{\tr}{trunk}
\usepackage{amsthm}

\usepackage{geometry}

\usepackage{fancyhdr}

\usepackage{color} 









\newcommand{\al}{\alpha}

\newcommand{\vphi}{\varphi}

\newcommand{\be}{\beta}

\newcommand{\ga}{\gamma}

\newcommand{\de}{\delta}

\newcommand{\om}{\omega}

\newcommand{\na}{\nabla}

\newcommand{\NA}{\nabla}

\newcommand{\bs}{\boldsymbol}

\newcommand{\ra}{\rightarrow}

\newcommand{\lra}{\longrightarrow}

\newcommand{\Ra}{\Rightarrow}

\newcommand{\xra}{\xrightarrow}

\newcommand{\xlra}{\xlongrightarrow}

\newcommand{\rgl}{\rangle}

\newcommand{\lgl}{\langle}

\newcommand{\dash}{\textrm{-}}

\newcommand{\ot}{\otimes}

\newcommand{\bpf}{\begin{proof}}

\newcommand{\epf}{\end{proof}}

\newcommand{\bthm}{\begin{thm}}

\newcommand{\ethm}{\end{thm}}

\newcommand{\bprop}{\begin{prop}}

\newcommand{\eprop}{\end{prop}}

\newcommand{\bcor}{\begin{cor}}

\newcommand{\ecor}{\end{cor}}

\newcommand{\blem}{\begin{lem}}

\newcommand{\elem}{\end{lem}}

\newcommand{\bdefn}{\begin{defn}}

\newcommand{\edefn}{\end{defn}}

\newcommand{\bexmp}{\begin{exmp}}

\newcommand{\eexmp}{\end{exmp}}

\newcommand{\brem}{\begin{rem}}

\newcommand{\erem}{\end{rem}}

\newcommand{\bdia}{\begin{displaymath}\xymatrix}

\newcommand{\edia}{\end{displaymath}}

\newcommand{\beq}{\begin{equation*}\begin{aligned}}

\newcommand{\eeq}{\end{aligned}\end{equation*}}

\newcommand{\bref}{\textbf{Ref}}

\newcommand{\intg}{\mathbb{Z}}

\newcommand{\real}{\mathbb{R}}

\newcommand{\comp}{\mathbb{C}}

\newcommand{\quot}{\mathbb{H}}

\newcommand{\afv}{\mathbb{A}}

\newcommand{\prv}{\mathbb{P}}

\newcommand{\mco}{\mathcal{O}}

\newcommand{\mcc}{\mathcal{C}}

\newcommand{\mcf}{\mathcal{F}}

\newcommand{\mcg}{\mathcal{G}}

\newcommand{\mcs}{\mathcal{S}}

\newcommand{\cp}{\mathbb{CP}}

\newcommand{\mfo}{\mathfrak{O}}

\newcommand{\mfg}{\mathfrak{g}}

\newcommand{\msa}{\mathscr{A}}

\newcommand{\msr}{\mathscr{R}}

\newcommand{\msg}{\mathscr{G}}

\newcommand{\msd}{\mathscr{D}}

\newcommand{\itbf}{\item\textbf}

\newcommand{\seqa}{a_1,...,a_}

\newcommand{\seqx}{x_1,...,x_}

\newcommand{\seqy}{y_1,...,y_}

\newcommand{\seqf}{f_1,...,f_}

\newcommand{\cred}{\textcolor{red}}

\newcommand{\cblue}{\textcolor{blue}}

\newcommand{\mfa}{\mathfrak{a}}

\newcommand{\mfb}{\mathfrak{b}}

\newcommand{\mfm}{\mathfrak{m}}

\newcommand{\mfn}{\mathfrak{n}}

\newcommand{\mfp}{\mathfrak{p}}

\newcommand{\Af}{A_{(f)}}


\newtheorem{thm}{\textbf {Theorem}}[section]

\newtheorem{cor}[thm]{\textbf{Corollary}}

\newtheorem{prop}[thm]{\textbf{Proposition}}

\newtheorem{lem}[thm]{\textbf{Lemma}}

\newtheorem{conj}[thm]{Conjecture}

\newtheorem{conv}[thm]{Convention}

\newtheorem{prob}[thm]{Problem}

\newtheorem{exer}[thm]{Exercise}

\newtheorem{quest}[thm]{Question}

\theoremstyle{definition}

\newtheorem{defn}[thm]{\textbf{Definition}}

\newtheorem{defns}[thm]{Definitions}

\newtheorem{exmp}[thm]{Example}

\newtheorem{exmps}[thm]{Examples}

\newtheorem{var}[thm]{Variant}

\newtheorem{vars}[thm]{Variants}

\newtheorem{con}[thm]{Construction}

\newtheorem{notn}[thm]{Notation}

\newtheorem{notns}[thm]{Notations}

\theoremstyle{remark}

\newtheorem{rem}[thm]{Remark}

\newtheorem{rems}[thm]{Remarks}

\newtheorem{warn}[thm]{Warning}

\newtheorem{sch}[thm]{Scholium}

\newtheorem{expl}[thm]{Explanations}

\newtheorem*{theorem}{\textbf{Theorem}}

\newtheorem*{corollary}{\textbf{Corollary}}

\newtheorem*{proposition}{\textbf{Proposition}}

\newtheorem*{lemma}{\textbf{Lemma}}

\newtheorem*{example}{\textbf{Example}}

\def\cok{\operatorname{Coker}}

\newcommand{\txi}{\tilde{\xi}}

\newcommand{\bxi}{\bar{\xi}}

\newcommand{\bz}{\bar{z}}

\begin{document}

\bibliographystyle{plain}

\else

\fi

\section{Balanced Pairings}

Recall from Definition \ref{balanced} that a \textit{balanced pairing} is a pairing $\mathcal{P}$ such that $\sigma(P) = 0.$ In this section, we prove the following theorem by induction:

\begin{thm}\label{balpair} For a fixed odd integer $n$, all balanced pairings of size $n$ are connected.
\end{thm}

To prove this theorem, we begin with the base case $n = 3$, which will be dealt with in Example \ref{exmp_base_case} below. Then, we use strong induction on $n,$ half of the number of vertices, assuming that all smaller odd values of $n$ satisfy Theorem \ref{balpair}. We divide the proof into two cases: planar and nonplanar pairings. After statements \ref{exmp_base_case} through \ref{manual}, we return to the proof of Theorem \ref{balpair} by using the pairing $\mathcal{P}_n$ to construct $D_{\mathcal{P}_n}.$ We then prove Theorem \ref{balpair} for our odd $n,$ completing the strong inductive argument.

\begin{exmp}\label{exmp_base_case}
Suppose $n = 3$. Then we have $6$ vertices on the boundary of the disk, so points of opposite parity can be paired in $3! = 6$ ways as below:

I. (1, 2), (3, 4), (5, 6). Signature = $2$.

II. (1, 2), (3, 6), (5, 4). Signature = $0$.

III. (1, 4), (3, 2), (5, 6). Signature = $0$.

IV. (1, 4), (3, 6), (5, 2). Signature = $0$.

V. (1, 6), (3, 2), (5, 4). Signature = $-2$.

VI. (1, 6), (3, 4), (5, 2). Signature = $0$.

Clearly there are four balanced pairings: II, III, IV, VI.

For pairing II, if we perform a cut-and-glue operation on the two pairs $(1,2)$ and $(4,5)$ then we get pairing IV. For pairing III if we perform a cut-and-glue operation on the two pairs $(2,3)$ and $(5,6)$, we also get pairing IV. For pairing VI, if we perform a cut-and-glue operation on the two pairs $(1, 6)$ and $(3, 4)$ we also get pairing IV. Thus all of the balanced pairings are connected to pairing IV, and therefore to each other.
\end{exmp}

\begin{exmp}\label{threelayers}
The argument that II, III, VI are all connected can be generalized. Suppose $(i,a),(j,b), (k,c)\in\mathcal{P}$ so that in the sequence $\{i,j,k,c,b,a\}$, adjacent indices will be of different parity. Then $\mathcal{P}$ is connected to two other pairings $\mathcal{P}'$ and $\mathcal{P}''$ where
$$(i,c),(j,k),(b,a)\in\mathcal{P}',\hspace{5mm}(i,j),(k,a),(b,c)\in\mathcal{P}'',$$
and the three pairings $\mathcal{P},\mathcal{P}', \mathcal{P}''$ coincide elsewhere.
\newline
The fact that $\mathcal{P}$ is connected to $\mathcal{P}'$ and $\mathcal{P}''$ occurs because we can do a cut-and-glue operation to get from one to the other. The coincidence elsewhere refers to the fact that the size of $\mathcal{P}$ is greater than $3,$ so the other vertices may match up. 

\end{exmp}
\begin{exmp}\label{onelayer} The cases I and V in Example \ref{exmp_base_case} can also be generalized as follows. Suppose a pairing $\mathcal{P}_n$ only has one layer. Then $|\sigma(\mathcal{P}_n)| = n - 1.$ In particular, $\mathcal{P}_n$ must be one of two pairings: either the unique pairing with a signature of $n - 1,$ or the unique pairing with a signature of $1 - n.$ In Figure \ref{onelayerfig} below, we see an example of this for $n = 3$:
\end{exmp}

\begin{figure}[H]
    \centering
    \includegraphics[scale = 0.15]{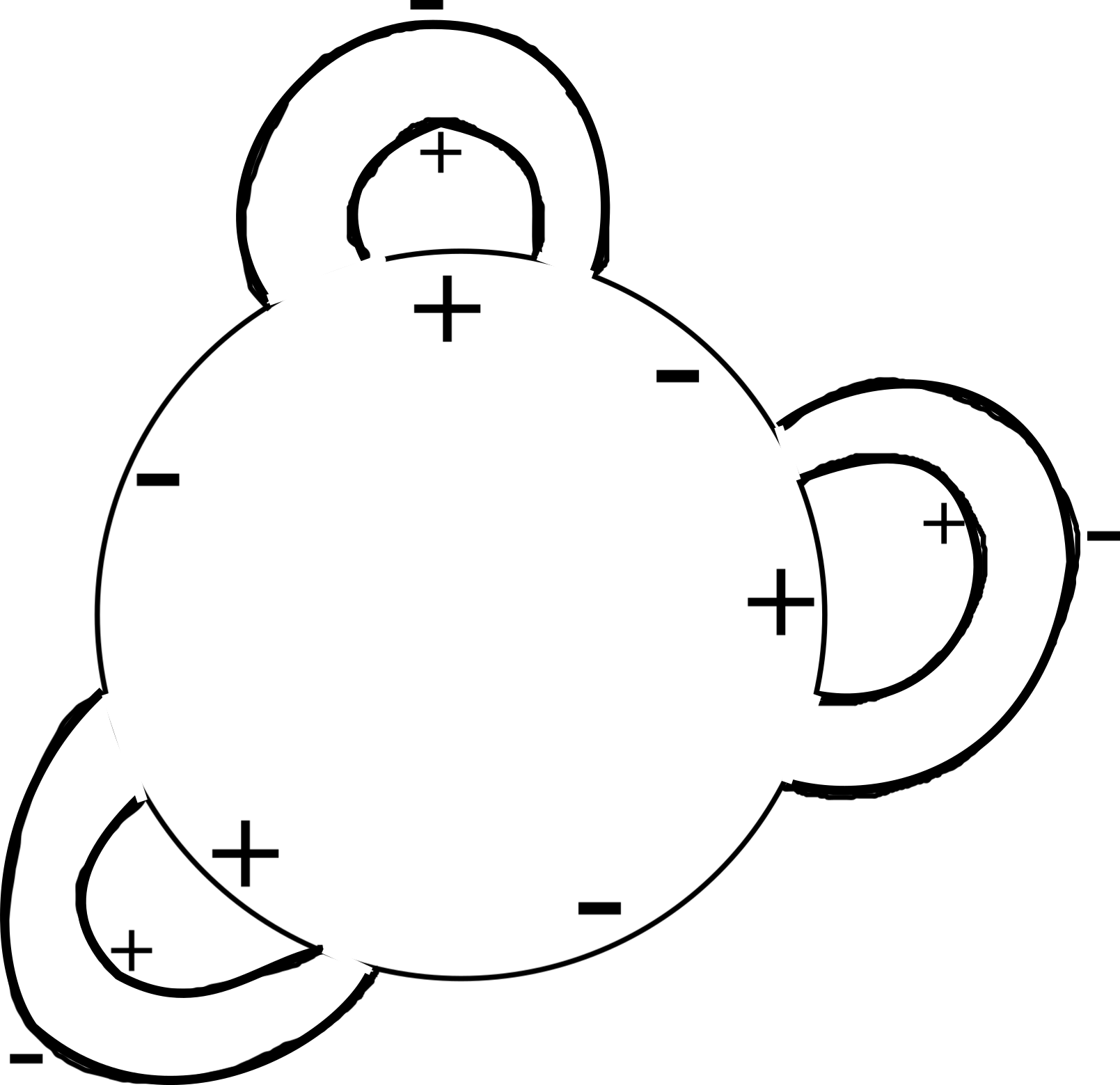}
    \caption{We have $\sigma(\mathcal{P}_3) = 2.$}
    \label{onelayerfig}
\end{figure}

Recall that $D$ is the unit disk on the complex plane $\mathbb{C}$ and with $2n$ vertices
$$p_k=e^{\frac{k-1}{n}2\pi i}.$$

Suppose $\mathcal{P}=\{(i_1,j_1),...,(i_n,j_n)\}$ is a planar pairing, then on $\mathbb{C}$ we can attach $n$ strips $S_1,...,S_n$ to $D$ so that $S_l$ is attached along $p_{i_l}$ and $p_{j_l}$. Pick an arc $\al_l \in S_l$ such that $\partial \al_l = \{P_{i_l}, P_{j_l} \}$ and we can orient $\al_l$ arbitrarily. We will call $\al_l$ {\it layers}. Consider $2n$ rays from the origin of $\mathbb{C}$: for $l=1,...,2n$, let
$$R_l=\{z|z=r\cdot e^{\frac{2l-1}{4n}2\pi i}~{\rm for}~{\rm some}~ r\in[0,\infty)\}.$$
Also give $R_l$ an arbitrary orientation.

\begin{figure}[H]
    \centering
    \includegraphics[scale = .15]{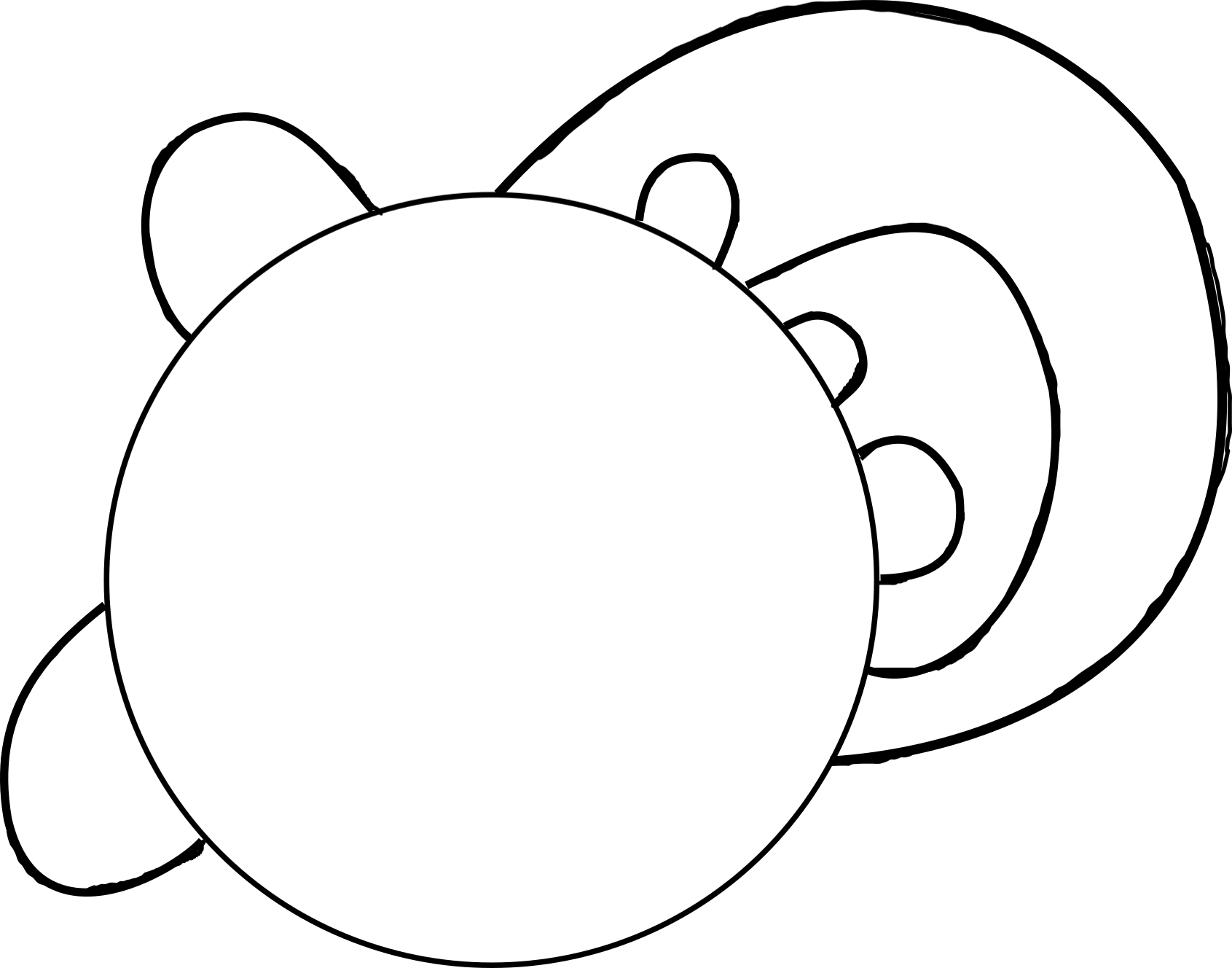}
    \caption{This pairing $\mathcal{P}_9$ has three layers. For simplicity, each arc denotes an entire strip.}
    \label{figlayers}
\end{figure}

In Figure \ref{figlayers}, we see visually why the pairing depicted has $3$ layers.

\begin{defn}\label{defn_layer}  
Under the above setting, for each $l$ with $1\leq l\leq 2n$, we define
$$x_l=\sum_{l=1}^{2n}|R_l\cdot \al_l|.$$
Here $|~|$ is the absolute value and $\cdot$ means the signed intersection of two curves.

Also we define
$$c(\mathcal{P}) = \sum_{l=1}^n x_l(\mathcal{P}).$$ 
Here $c(\mathcal{P})$ measures to some extent the complexity of the pairing $\mathcal{P}$. 

\end{defn}

\brem
The definition of $x_l$ and $c$ does not only depend on the balanced pairing but also how $D_{\mathcal{P}}$ is embedded into $\mathbb{C}$. However, the different ways of embedding do not make a difference in our argument.
\erem

\begin{lem}\label{redcomp} Suppose we have a pairing 
$$\mathcal{P}=\{(i_1,j_1),...,(i_n,j_n)\}.$$
If some $x_l(\mathcal{P})\geq 3$, then we can apply two cut-and-glue operations (as in Example \ref{threelayers}) to get a new pairing $\mathcal{P}'$ where $c(\mathcal{P}') < c(\mathcal{P})$ and $x_l(\mathcal{P'}) < x_l(\mathcal{P}).$

\end{lem}

\begin{proof} Consider the outermost three layers of $D_\mathcal{P}$ that intersect $R_l.$ Without loss of generality, we can assume that they are $\al_1$, $\al_2$ and $\al_3$. Note $i_1$ and $i_2$ must be of different parity or there must be an odd number of vertices in between $\al_1$ and $\al_2$. If this is the case, then at least one of the vertices must connect to another vertex to form the second outermost layer, contradicting the assumption that we begin by considering the outermost three layers of $D_\mathcal{P}.$ Similarly $i_2$ and $i_3$ must be of different parity. So we can apply two cut-and-glue operations as in Example \ref{threelayers} and it is straightforward to check that the new pairing $\mathcal{P}'$ satisfies the requirement of the lemma.

\end{proof}

\begin{cor}\label{reducecomplexity} Given any planar balanced pairing $\mathcal{P}$, we can find another planar balanced pairing $\mathcal{P}'$, so that $\mathcal{P}$ and $\mathcal{P}'$ are connected and we have $x_l(\mathcal{P}') \leq 2$ for all $l.$
\end{cor}

\begin{proof}
It is a fact that $c(\mathcal{P})>n$ because there are $n$ strips in total and each strip increases the sum of the $x_l$ by at least $1$.

By Lemma \ref{redcomp}, we can decrease $c(\mathcal{P})$ whenever $x_{l}(\mathcal{P})\geq 3$ for some $l$. However, we must always have $c(\mathcal{P}) \geq n$ so we cannot decrease it forever, and thus this decreasing process will stop when $x_{l}(\mathcal{P}')\leq 2$ for all $l$.
\end{proof}

For the following lemma, we recall Definition \ref{pairingdefn}, where we defined $\Pi_n.$ 

\begin{lem}\label{ideal} Suppose a balanced pairing $\mathcal{P}_n$ has two pairs $(i, i+ 1)$ and $(n + i, n + i + 1)$. Then $\mathcal{P}_n$ is connected to all other balanced pairings of size $n$ which also contain these two pairs.

\end{lem}

\begin{figure}[H]
    \centering
    \includegraphics[scale = .15]{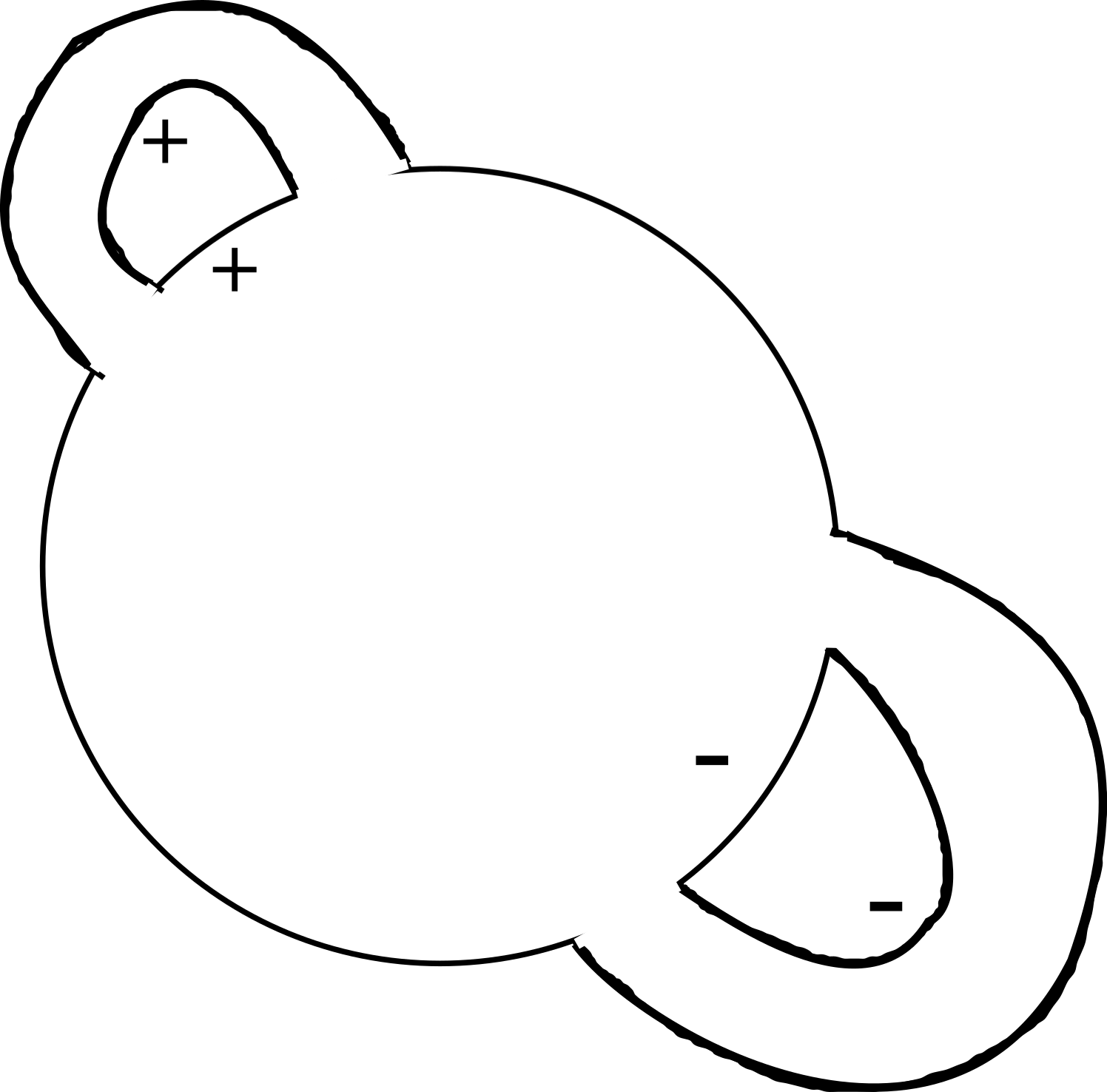}
    \caption{The basic idea is that we remove the positive and negative boundary components whenever we have a pair as in the figure above, allowing us to use induction.}
    \label{posnegfig}
\end{figure}

\bpf We form a new pairing $\mathcal{P}_{n - 2}$ by removing the pairs $(i, i + 1), (n + i, n + i + 1)$ from $\mathcal{P}_n$ and leaving all other pairs unchanged. Since $\sigma(\mathcal{P}_n) = \sigma(\mathcal{P}_{n - 2}) = 0,$ we know from the inductive hypothesis that $\mathcal{P}_{n - 2}$ is connected to all other balanced pairings of size $n - 2.$ There is a bijection between each of the balanced pairings of size $n - 2$ and balanced pairings of size $n$ with the pairs $(i, i + 1), (n + i, n + i + 1).$ 

After removing these two pairs, we form $P_s$ of size $n - 2$ and then add the two pairs back. After we do a cut-and-glue operation on the pairs $(i, i + 1), (n + i, n + i + 1)$ we get $\mathcal{P}_s$ of size $n$, as defined in Lemma \ref{special}. Therefore, all balanced pairings of size $n$ that contain pairs of the form $(i, i + 1)$ and $(n + i, n + i + 1)$ are connected to $\mathcal{P}_s$ of size $n,$ so they are all connected to each other. 



\epf


 




\begin{cor}\label{idealinduction} Suppose we have a pairing $\mathcal{P}_n \in \Pi_n$ which has two pairs $(i, i + 1), (j, j + 1)$ where $i, j$ have different parity. Then this pairing is connected to $\mathcal{P}_s$ of size $n.$
\end{cor}

\bpf

Note that $\mathcal{P}_s$ is rotation invariant so we can relabel the vertices so that $i=k+1$ where $n = 2k + 1$. If $j = n + k + 1,$ then we are done by Lemma \ref{ideal}.

If $j\neq k+n+1$, then the assumption is that $i, j$ are not the same parity so $j$ and $k+n+1$ are in the same parity as $n$ is odd. Since we have the freedom to choose any balanced pairing $\mathcal{P}_{n - 2}$ that we would like, we can choose one that has a strip which connects $n + k + 1$ to $n + k + 2.$ We know that the inner boundary component created by this strip will have the same sign as that created by the strip connecting $p_j$ to $p_{j + 1}$, and therefore the opposite sign from the boundary component created by the strip connecting $p_i$ and $p_{i + 1}.$ Then, we can replace $j$ by $k + n + 1$ and we are done by Lemma \ref{ideal}.
\epf

These allow us to prove the following theorem:

\begin{thm}\label{planpair} For a fixed $n,$ all balanced planar pairings of size $n$ are connected.
\end{thm}

\begin{proof}
Suppose $\mathcal{P}_n\in\Pi_n$ is any planar balanced pairing of size $n$.

Case 1: There exist integers $i, j$ where $1 \leq i < j \leq 2n$ which satisfy the conditions of Corollary \ref{idealinduction}. Then we are done by that corollary.

Case 2: If $(i, i + 1), (j, j + 1) \in \mathcal{P}_n$, then $i \equiv j \pmod 2.$

Corollary \ref{reducecomplexity} implies that we can assume $x_l(\mathcal{P}_n) \leq 2$ for all $l$ satisfying $1 \leq l \leq 2n.$ 
Without loss of generality, we can assume that $(1,2p)\in \mathcal{P}_n$, where $p\neq 1,n$. (If $n>1$ and all of the pairs in $\mathcal{P}_n$ are of the form $(i,i+1)$, then $\mathcal{P}_n$ cannot be balanced, see Example \ref{onelayer}.)

Then we claim the following things:

(1). We have $(2p-2,2p-1)\in \mathcal{P}_n.$

(2). There exists a $q$ so that $q\neq p+1$ and $(2p + 1, 2q) \in \mathcal{P}_n$.

(3). We have $(2p + 2, 2p + 3) \in \mathcal{P}_n.$

To prove the claim, (1) follows from the fact that $x_{2p-2}(\mathcal{P}_n)\leq 2$, where $x_{2p-2}$ is defined in Definition \ref{defn_layer}. If (2) is not true, then $(2p+1,2p+2)\in \mathcal{P}_n$, but then $2p-2\not\equiv 2p+1 \mod 2$ and this violates the assumption of case 2. Similarly (3) follows from the fact that $x_{2p+2}(\mathcal{P})\leq 2$.

Then we apply two cut and paste as in Example \ref{threelayers}. We know that $\mathcal{P}_{n}$ is connected to another planar balanced pairing $\mathcal{P}'_n$ where $(1,2q),(2p-2,2p+1),(2p+2,2p+3)\in\mathcal{P}'_n$. Then we can apply Corollary \ref{idealinduction} and finish the proof. 




\end{proof}




\brem While so far we have only been working on a disk, we never used any property of the disk other than that it is a surface with exactly one boundary. Therefore, all of our above results for disks are also true for any surface with connected boundary.

\erem



\begin{lem}\label{manual} All balanced pairings of size $n,$ including nonplanar ones, are connected for $n = 5, 7.$

\begin{proof} This follows from directly checking all $68$ balanced pairings when $n = 5$ and all $2588$ balanced pairings when $n = 7$ with a computer program.
 
\end{proof}

\end{lem}

Now we are ready to prove Theorem \ref{balpair}, which states that all balanced pairings, including nonplanar ones, for a fixed $n$ are connected. 
\begin{proof}

We know that all balanced, planar pairings are connected from Theorem \ref{planpair}. We claim that all balanced nonplanar pairings can be transformed into planar pairings and will prove this by strong induction. The base case $n=3$ is addressed in Example \ref{exmp_base_case}. For $n>3$ and $n$ odd, we proceed as follows.


Suppose we have a general nonplanar pairing $\mathcal{P}_n$. By assumption, since the pairing is nonplanar, there exist two strips $S_i$ and $S_j$ on the circle which intersect each other in the geometric embedding into $\mathbb{C}$. Then the union of the disk and these two strips forms a surface with connected boundary as shown in Figure \ref{oneboundfig}. On this connected boundary, the pairing $\mathcal{P}_n$ naturally induces a pairing of size $n - 2$ once we remove the four vertices connected by $S_i, S_j$.

\begin{figure}[H]
    \centering
    \includegraphics[scale = 0.15]{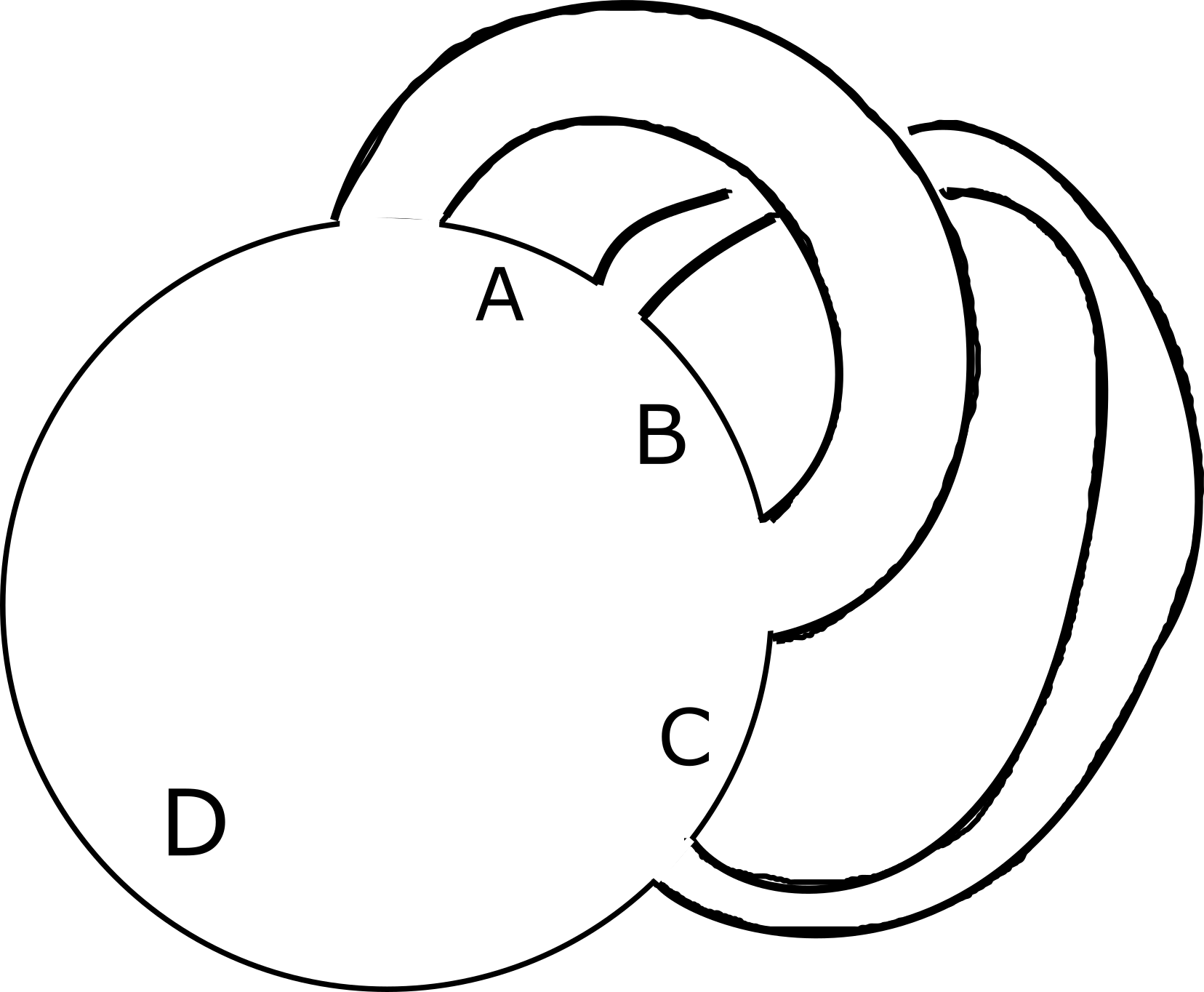}
    \caption{This is a surface with one boundary.}
    \label{oneboundfig}
\end{figure}

As shown in Figure \ref{oneboundfig}, we have four regions A, B, C, D that are separated by the vertices which are the boundaries of the two strips $S_i, S_j.$ Suppose one of these regions has 5 vertices. Then, since we have the freedom to choose any balanced pairing $\mathcal{P}_{n - 2},$ we pick one where two arcs as in Corollary \ref{idealinduction} are in the region with $5$ vertices. 

Now consider if there are $4$ vertices in region A and 2 vertices in regions B, C, D. Suppose we then add one more vertex. If we add it to region $A,$ we have $5$ vertices in one region, so we are done by the argument above. Otherwise, we can choose one of the pairs in A along with one of the pairs in either B or D and apply Corollary \ref{idealinduction}. 

For this to happen, there must be $4 + 4 + 2 + 2 + 4 + 1 = 17$ vertices total, implying that $n \geq 9.$ Combining this with Example \ref{exmp_base_case} and Lemma \ref{manual}, we get the result for all odd integers $n.$

\end{proof}

\ifx\allfiles\undefined

\bibliography{Index}

\end{document}

\fi

\ifx\allfiles\undefined

\documentclass[12pt,a4paper]{article}


\usepackage{graphicx}

\usepackage{amsmath}

\usepackage{amssymb}
\DeclareMathOperator{\tr}{trunk}
\usepackage{amsthm}

\usepackage{geometry}

\usepackage{fancyhdr}

\usepackage{color} 









\newcommand{\al}{\alpha}

\newcommand{\vphi}{\varphi}

\newcommand{\be}{\beta}

\newcommand{\ga}{\gamma}

\newcommand{\de}{\delta}

\newcommand{\om}{\omega}

\newcommand{\na}{\nabla}

\newcommand{\NA}{\nabla}

\newcommand{\bs}{\boldsymbol}

\newcommand{\ra}{\rightarrow}

\newcommand{\lra}{\longrightarrow}

\newcommand{\Ra}{\Rightarrow}

\newcommand{\xra}{\xrightarrow}

\newcommand{\xlra}{\xlongrightarrow}

\newcommand{\rgl}{\rangle}

\newcommand{\lgl}{\langle}

\newcommand{\dash}{\textrm{-}}

\newcommand{\ot}{\otimes}

\newcommand{\bpf}{\begin{proof}}

\newcommand{\epf}{\end{proof}}

\newcommand{\bthm}{\begin{thm}}

\newcommand{\ethm}{\end{thm}}

\newcommand{\bprop}{\begin{prop}}

\newcommand{\eprop}{\end{prop}}

\newcommand{\bcor}{\begin{cor}}

\newcommand{\ecor}{\end{cor}}

\newcommand{\blem}{\begin{lem}}

\newcommand{\elem}{\end{lem}}

\newcommand{\bdefn}{\begin{defn}}

\newcommand{\edefn}{\end{defn}}

\newcommand{\bexmp}{\begin{exmp}}

\newcommand{\eexmp}{\end{exmp}}

\newcommand{\brem}{\begin{rem}}

\newcommand{\erem}{\end{rem}}

\newcommand{\bdia}{\begin{displaymath}\xymatrix}

\newcommand{\edia}{\end{displaymath}}

\newcommand{\beq}{\begin{equation*}\begin{aligned}}

\newcommand{\eeq}{\end{aligned}\end{equation*}}

\newcommand{\bref}{\textbf{Ref}}

\newcommand{\intg}{\mathbb{Z}}

\newcommand{\real}{\mathbb{R}}

\newcommand{\comp}{\mathbb{C}}

\newcommand{\quot}{\mathbb{H}}

\newcommand{\afv}{\mathbb{A}}

\newcommand{\prv}{\mathbb{P}}

\newcommand{\mco}{\mathcal{O}}

\newcommand{\mcc}{\mathcal{C}}

\newcommand{\mcf}{\mathcal{F}}

\newcommand{\mcg}{\mathcal{G}}

\newcommand{\mcs}{\mathcal{S}}

\newcommand{\cp}{\mathbb{CP}}

\newcommand{\mfo}{\mathfrak{O}}

\newcommand{\mfg}{\mathfrak{g}}

\newcommand{\msa}{\mathscr{A}}

\newcommand{\msr}{\mathscr{R}}

\newcommand{\msg}{\mathscr{G}}

\newcommand{\msd}{\mathscr{D}}

\newcommand{\itbf}{\item\textbf}

\newcommand{\seqa}{a_1,...,a_}

\newcommand{\seqx}{x_1,...,x_}

\newcommand{\seqy}{y_1,...,y_}

\newcommand{\seqf}{f_1,...,f_}

\newcommand{\cred}{\textcolor{red}}

\newcommand{\cblue}{\textcolor{blue}}

\newcommand{\mfa}{\mathfrak{a}}

\newcommand{\mfb}{\mathfrak{b}}

\newcommand{\mfm}{\mathfrak{m}}

\newcommand{\mfn}{\mathfrak{n}}

\newcommand{\mfp}{\mathfrak{p}}

\newcommand{\Af}{A_{(f)}}


\newtheorem{thm}{\textbf {Theorem}}[section]

\newtheorem{cor}[thm]{\textbf{Corollary}}

\newtheorem{prop}[thm]{\textbf{Proposition}}

\newtheorem{lem}[thm]{\textbf{Lemma}}

\newtheorem{conj}[thm]{Conjecture}

\newtheorem{conv}[thm]{Convention}

\newtheorem{prob}[thm]{Problem}

\newtheorem{exer}[thm]{Exercise}

\newtheorem{quest}[thm]{Question}

\theoremstyle{definition}

\newtheorem{defn}[thm]{\textbf{Definition}}

\newtheorem{defns}[thm]{Definitions}

\newtheorem{exmp}[thm]{Example}

\newtheorem{exmps}[thm]{Examples}

\newtheorem{var}[thm]{Variant}

\newtheorem{vars}[thm]{Variants}

\newtheorem{con}[thm]{Construction}

\newtheorem{notn}[thm]{Notation}

\newtheorem{notns}[thm]{Notations}

\theoremstyle{remark}

\newtheorem{rem}[thm]{Remark}

\newtheorem{rems}[thm]{Remarks}

\newtheorem{warn}[thm]{Warning}

\newtheorem{sch}[thm]{Scholium}

\newtheorem{expl}[thm]{Explanations}

\newtheorem*{theorem}{\textbf{Theorem}}

\newtheorem*{corollary}{\textbf{Corollary}}

\newtheorem*{proposition}{\textbf{Proposition}}

\newtheorem*{lemma}{\textbf{Lemma}}

\newtheorem*{example}{\textbf{Example}}

\def\cok{\operatorname{Coker}}

\newcommand{\txi}{\tilde{\xi}}

\newcommand{\bxi}{\bar{\xi}}

\newcommand{\bz}{\bar{z}}

\begin{document}

\bibliographystyle{plain}

\else

\fi

\section{Surfaces with Multiple Boundary Components}

We now extend the results of Section 4 to surfaces with more than one boundary component. As before, we require that each boundary component has an even number of vertices and the endpoints of any strip are vertices whose indices are of opposite parity. Unlike in Section 4, we can now join vertices with strips that are on separate boundary components.

We first deal with $2$ boundary components, and then we generalize to any number of boundary components through a strong inductive argument and by noting that the strips can be thought of as part of the surface as we saw in Figure \ref{oneboundfig}.

\begin{notn}\label{notn1} We use $S$ to denote our surface with multiple boundary components, and we say that $\partial S = C_1 \cup C_2 \cup \ldots \cup C_s,$ where the $C_i$ denote the boundary components of $S.$ 

\end{notn}

\begin{notn} Analogously to Section 3, we give $S$ an orientation which induces an orientation on $\partial S$ as well. Each boundary component $C_i$ defined in Notation \ref{notn1} has vertices $p_{2n_{i - 1} + 1}, \ldots p_{2n_i}.$ For all of these boundary components, we assign the minor arc from a vertex of odd index to a vertex of even index (according to the orientation of $\partial{S}$) to be positive.
\end{notn}

\begin{rem}
We require that there are an even number of vertices on each boundary components and this actually follows from the fact that if $(M,\ga)$ is a balanced sutured manifold as introduced in the introduction and $S$ is a properly embedded surfaces in $M$, then any boundary component of $S$ intersects the suture $\ga$ an even number of times.
\end{rem}

\begin{lem}\label{sameparity} For a balanced pairing to exist, we must have $s + n_1 + \ldots + n_s \equiv 0 \pmod 2.$ 

\end{lem}

\bpf The original surface $S$ has $s$ boundary components and $n_1 + n_2 + \ldots + n_s$ pairs of vertices, where each pair is connected by exactly one strip. Then $\chi(S) = 2 - 2g(S) - s.$ 

After we attach all of the strips, we get the new surface $S'$ which has $b$ boundary components. In order to have a balanced pairing, we require $b$ to be even. Therefore, we have $\chi(S') = 2 - 2g(S') - b = 2 - 2g(S) - (s +n_1 + n_2 + \ldots + n_s).$ It follows that $s + n_1 + \ldots + n_s \equiv 0 \pmod 2.$

\epf

\begin{lem}\label{even} Suppose we have a surface $S$ with $2$ boundary components. Then the number of strips between the $2$ boundary components is even. Further, for each of the boundary components, the number of vertices with odd index and the number of vertices with even index connected by a strip to the other boundary component are the same.

\bpf Suppose by contradiction that there are an odd number of strips between the two boundary components. Then, that leaves an odd number of vertices on each of the boundary components to be joined in pairs by strips, which is impossible.

Now suppose that the number of vertices with odd index and the number of vertices with even index connected by a strip to the other boundary component are not the same. Then consider the strips connecting the vertices on that specific boundary component. Since the number of vertices of odd and even index are different, it is impossible to pair all of the vertices with strips that connect vertices with odd indices to vertices with even indices.

\epf

\end{lem}

\begin{lem}\label{twoneeded} For a surface $S$ with exactly two boundary components, any two balanced pairings are connected.

\bpf Suppose $\mathcal{P}_n$ is a balanced pairing of size $n=n_1+n_2$. Without loss of generality, we let $n_2 = n_1 + 2k,$ where $k \geq 0.$ We prove in three steps that $\mathcal{P}_n$ can be connected to a modeled balanced pairing and hence any two balanced pairings are connected through the modeled one.

{\bf Step 1.} If $C_1, C_2$ are separate boundary components without any strips between them, we can do a cut-and-glue operation on one strip on $C_1$ and another strip on $C_2$ to get two strips between $C_1$ and $C_2.$ Then we can relabel the vertices on $C_1, C_2$ to assume that $p_1$ is connected to $p_{2n_1 + 2}.$ We call the strip connecting these two vertices $S_1.$

{\bf Step 2.} Let $S' = S \cup S_1.$ Then $S'$ has connected boundary as shown in Figure \ref{fig9}:

\begin{figure}[H]
    \centering
    \includegraphics[scale = .15]{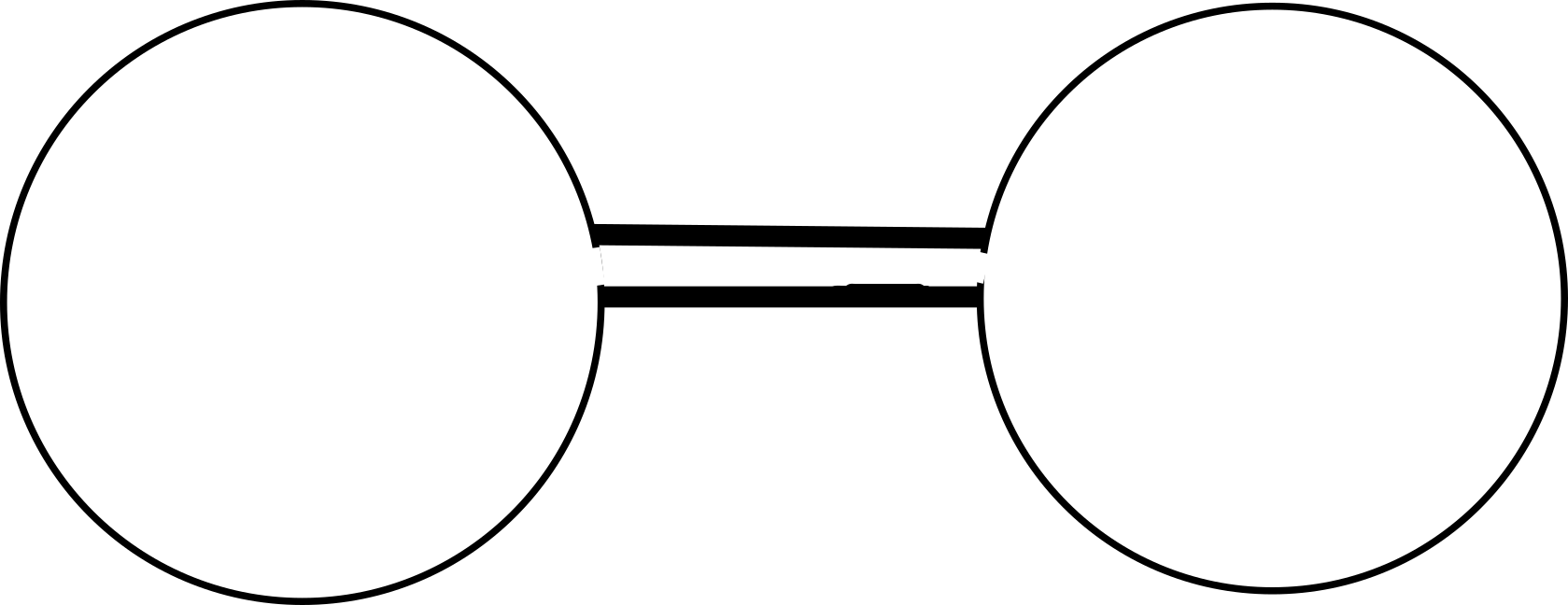}
    \caption{The boundary of $S'$ is connected.}
    \label{fig9}
\end{figure}

Therefore, we can apply Theorem \ref{balpair} to get that all balanced pairings on $S'$ are connected. For this reason, we can say that our original pairing $\mathcal{P}_n$ is connected to a new pairing $\mathcal{P}'_n$ where the following holds:

1. The vertex $p_1$ is connected to $p_{2n_1 + 2}$ with the strip $S_1$ and the vertex $p_2$ is connected to $p_{4n_1 + 2k + 1}$ with the strip $S_2.$ 

2. All other vertices $p_i$ are connected to $p_{i + 1}.$

In Figure \ref{redblueblack}, we show an example of $\mathcal{P}'_n,$ which in general satisfies the same addition relation as in this figure.

\begin{figure}[H]
    \centering
    \includegraphics[scale = .15]{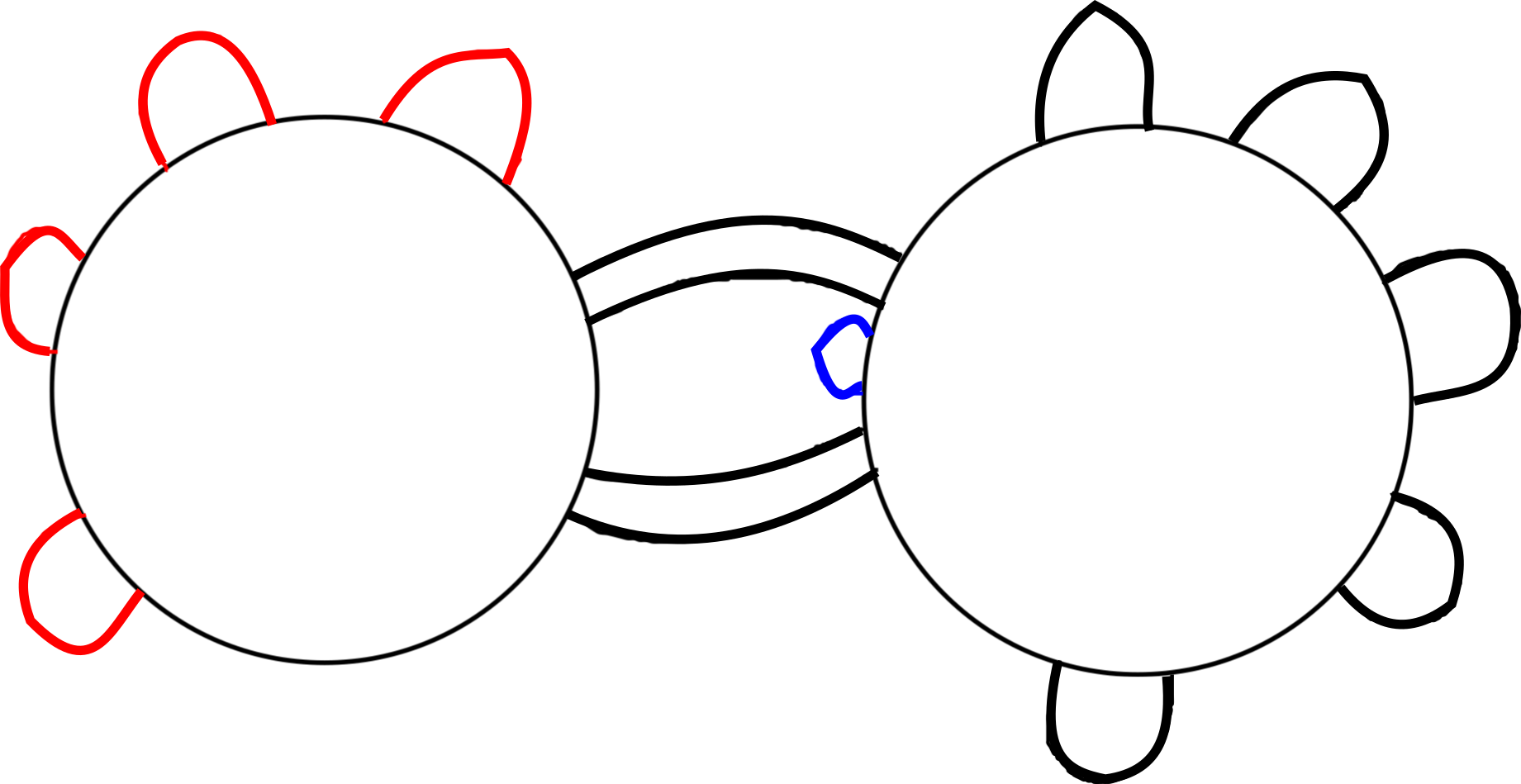}
    \caption{Note that there are $4$ red strips on the left circle, $1$ blue strip on the left side of the right circle, and $5$ black strips on the right side of the right circle. We have $4 + 1 = 5.$}
    \label{redblueblack}
\end{figure}

By a direct computation, we have that $\mathcal{P}'_n$ is balanced.

Then, we can perform a cut-and-glue operation along $S_1$ and $S_2$ to separate $C_1$ and $C_2.$ Thus $\mathcal{P}_n$ is connected to $\mathcal{P}_n''=\mathcal{P}_{n_1}\sqcup\mathcal{P}_{n_2}$, which means that $\mathcal{P}_n''$ splits into two separate parings $\mathcal{P}_{n_1}$ and $\mathcal{P}_{n_2}$ of size $n_1$ and $n_2$ respectively. 

{\bf Step 3.} Note that $\sigma(\mathcal{P}_{n_1}) = n - 1$, so $\mathcal{P}_{n - 1}$ is unique by the same argument as in Example \ref{onelayer}. However, $\mathcal{P}_{n_2}$ is not unique because of the other possibilities that exist due to the relabeling in Step 1. To address this, we use two cut-and-glue operations as in Example \ref{exmp_base_case}, which results in a rotation of $\mathcal{P}_{n_2}.$ This exactly resolves the ambiguity due to the relabeling. So we are done.

\epf

\end{lem}

\begin{thm}\label{balpairall} For a fixed surface $S$, all balanced pairings are connected.

\end{thm}

\bpf
We prove this theorem by strong induction on the number of boundary components of $S$. Say
$$\partial{S}=C_1\cup...\cup C_s.$$

We have already addressed the base cases where $s = 1$ and $s = 2.$ Now we will assume the hypothesis holds for all positive integers up to $s$ and prove it for $s + 1.$ Suppose $\mathcal{P}_n$ is a balanced pairing of size $n$.

Note each boundary component $C_i$ has $2n_i$ vertices. We have two cases: either at least one of the $n_i$ is odd or all of the $n_i$ are even.

{\bf Case 1}: At least one of the $n_i$ is odd. Without loss of generality, let us assume that $n_1=2k+1.$

As in the proof of Lemma \ref{twoneeded}, we can assume that $C_1$ is connected to at least one of the other boundary components. After relabeling, we can assume that the vertex $p_1$ is connected to $p_{2n_1+2}$ on $C_2$ via a strip $S_1$. As in the proof of Lemma \ref{twoneeded}, we can look at $S\cup S_1$ which has $s-1$ many boundary components and the inductive hypothesis applies. As a result we can assume that $\mathcal{P}_n$ is connected to $\mathcal{P}_n'$ where on $C_1$, $\mathcal{P}_n'$ has the following form:

(1). $p_i$ is connected to $p_{i+1}$ unless $i=1$ or $i=2k+2$.

(2). $p_1$ is connected to $p_{2n_1+2}$ via the strip $S_1$ and $p_{2k+2}$ is connected to $p_{2n_1+1}$ via the strip $S_2$.

Then, we can perform a cut-and-glue operation along $S_1$ and $S_2$ to isolate $C_i$ from the other boundary components. Then the strips connected to $C_i$ are of the form depicted below:

\begin{figure}[H]
    \centering
    \includegraphics[scale = .15]{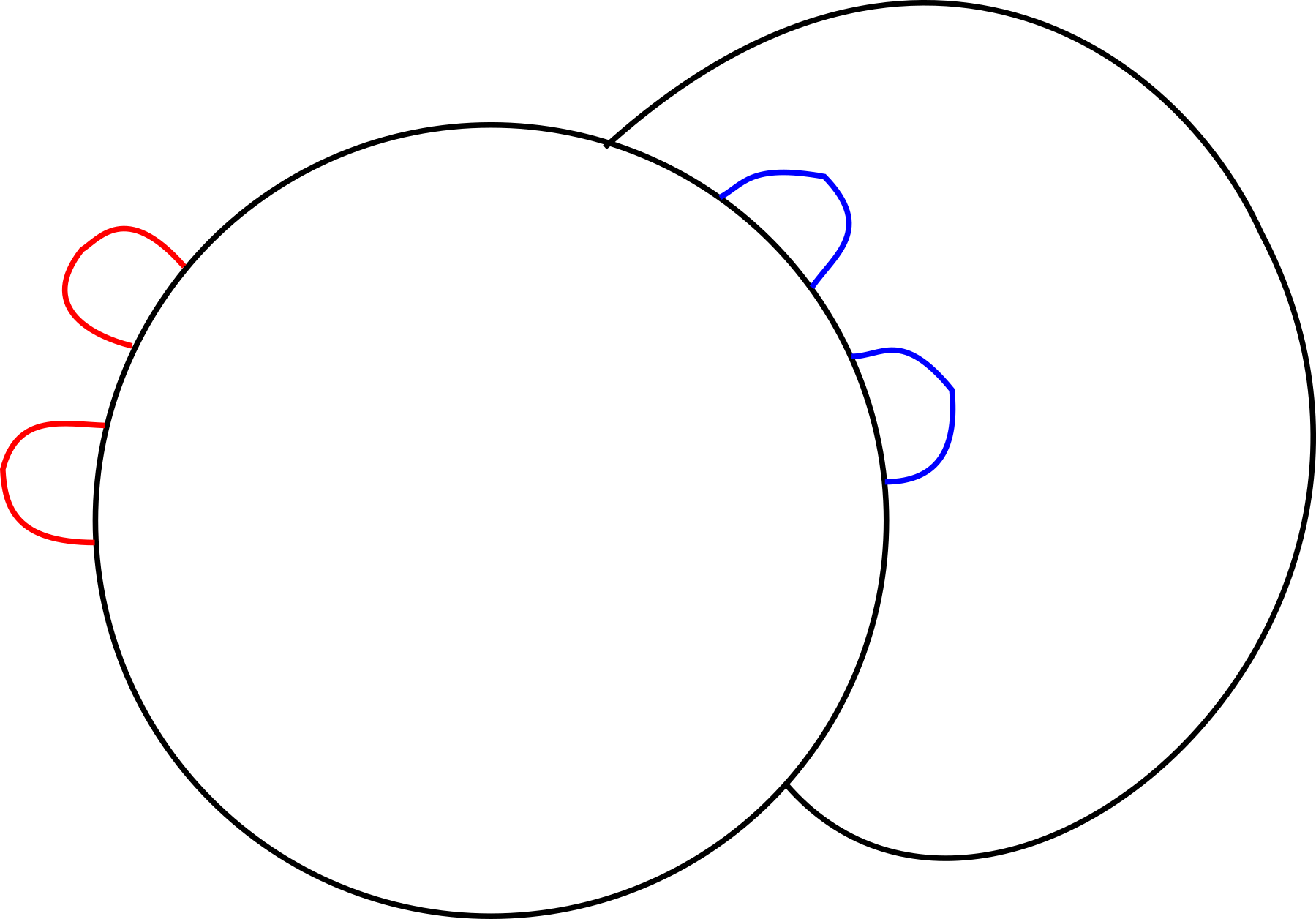}
    \caption{The number of red strips on the left side equals the number of blue strips on the right side.}
    \label{fig:my_label}
\end{figure}

Therefore, $\sigma(C_i) = 0,$ so we are done by strong induction.

{\bf Case 2}: All of the $n_i$ are even.

As before, we can view two boundary components connected by a strip as a surface with connected boundary. Then the number of vertices on this new boundary component is congruent to $2 \mod 4,$ so we are done by Case 1 and Lemma \ref{twoneeded}.

\epf

This gives us Proposition \ref{multbound}, which proves Corollary \ref{topapp}.

\ifx\allfiles\undefined

\bibliography{Index}

\end{document}

\fi

\ifx\allfiles\undefined

\documentclass[12pt,a4paper]{article}


\usepackage{graphicx}

\usepackage{amsmath}

\usepackage{amssymb}
\DeclareMathOperator{\tr}{trunk}
\usepackage{amsthm}

\usepackage{geometry}

\usepackage{fancyhdr}

\usepackage{color} 









\newcommand{\al}{\alpha}

\newcommand{\vphi}{\varphi}

\newcommand{\be}{\beta}

\newcommand{\ga}{\gamma}

\newcommand{\de}{\delta}

\newcommand{\om}{\omega}

\newcommand{\na}{\nabla}

\newcommand{\NA}{\nabla}

\newcommand{\bs}{\boldsymbol}

\newcommand{\ra}{\rightarrow}

\newcommand{\lra}{\longrightarrow}

\newcommand{\Ra}{\Rightarrow}

\newcommand{\xra}{\xrightarrow}

\newcommand{\xlra}{\xlongrightarrow}

\newcommand{\rgl}{\rangle}

\newcommand{\lgl}{\langle}

\newcommand{\dash}{\textrm{-}}

\newcommand{\ot}{\otimes}

\newcommand{\bpf}{\begin{proof}}

\newcommand{\epf}{\end{proof}}

\newcommand{\bthm}{\begin{thm}}

\newcommand{\ethm}{\end{thm}}

\newcommand{\bprop}{\begin{prop}}

\newcommand{\eprop}{\end{prop}}

\newcommand{\bcor}{\begin{cor}}

\newcommand{\ecor}{\end{cor}}

\newcommand{\blem}{\begin{lem}}

\newcommand{\elem}{\end{lem}}

\newcommand{\bdefn}{\begin{defn}}

\newcommand{\edefn}{\end{defn}}

\newcommand{\bexmp}{\begin{exmp}}

\newcommand{\eexmp}{\end{exmp}}

\newcommand{\brem}{\begin{rem}}

\newcommand{\erem}{\end{rem}}

\newcommand{\bdia}{\begin{displaymath}\xymatrix}

\newcommand{\edia}{\end{displaymath}}

\newcommand{\beq}{\begin{equation*}\begin{aligned}}

\newcommand{\eeq}{\end{aligned}\end{equation*}}

\newcommand{\bref}{\textbf{Ref}}

\newcommand{\intg}{\mathbb{Z}}

\newcommand{\real}{\mathbb{R}}

\newcommand{\comp}{\mathbb{C}}

\newcommand{\quot}{\mathbb{H}}

\newcommand{\afv}{\mathbb{A}}

\newcommand{\prv}{\mathbb{P}}

\newcommand{\mco}{\mathcal{O}}

\newcommand{\mcc}{\mathcal{C}}

\newcommand{\mcf}{\mathcal{F}}

\newcommand{\mcg}{\mathcal{G}}

\newcommand{\mcs}{\mathcal{S}}

\newcommand{\cp}{\mathbb{CP}}

\newcommand{\mfo}{\mathfrak{O}}

\newcommand{\mfg}{\mathfrak{g}}

\newcommand{\msa}{\mathscr{A}}

\newcommand{\msr}{\mathscr{R}}

\newcommand{\msg}{\mathscr{G}}

\newcommand{\msd}{\mathscr{D}}

\newcommand{\itbf}{\item\textbf}

\newcommand{\seqa}{a_1,...,a_}

\newcommand{\seqx}{x_1,...,x_}

\newcommand{\seqy}{y_1,...,y_}

\newcommand{\seqf}{f_1,...,f_}

\newcommand{\cred}{\textcolor{red}}

\newcommand{\cblue}{\textcolor{blue}}

\newcommand{\mfa}{\mathfrak{a}}

\newcommand{\mfb}{\mathfrak{b}}

\newcommand{\mfm}{\mathfrak{m}}

\newcommand{\mfn}{\mathfrak{n}}

\newcommand{\mfp}{\mathfrak{p}}

\newcommand{\Af}{A_{(f)}}


\newtheorem{thm}{\textbf {Theorem}}[section]

\newtheorem{cor}[thm]{\textbf{Corollary}}

\newtheorem{prop}[thm]{\textbf{Proposition}}

\newtheorem{lem}[thm]{\textbf{Lemma}}

\newtheorem{conj}[thm]{Conjecture}

\newtheorem{conv}[thm]{Convention}

\newtheorem{prob}[thm]{Problem}

\newtheorem{exer}[thm]{Exercise}

\newtheorem{quest}[thm]{Question}

\theoremstyle{definition}

\newtheorem{defn}[thm]{\textbf{Definition}}

\newtheorem{defns}[thm]{Definitions}

\newtheorem{exmp}[thm]{Example}

\newtheorem{exmps}[thm]{Examples}

\newtheorem{var}[thm]{Variant}

\newtheorem{vars}[thm]{Variants}

\newtheorem{con}[thm]{Construction}

\newtheorem{notn}[thm]{Notation}

\newtheorem{notns}[thm]{Notations}

\theoremstyle{remark}

\newtheorem{rem}[thm]{Remark}

\newtheorem{rems}[thm]{Remarks}

\newtheorem{warn}[thm]{Warning}

\newtheorem{sch}[thm]{Scholium}

\newtheorem{expl}[thm]{Explanations}

\newtheorem*{theorem}{\textbf{Theorem}}

\newtheorem*{corollary}{\textbf{Corollary}}

\newtheorem*{proposition}{\textbf{Proposition}}

\newtheorem*{lemma}{\textbf{Lemma}}

\newtheorem*{example}{\textbf{Example}}

\def\cok{\operatorname{Coker}}

\newcommand{\txi}{\tilde{\xi}}

\newcommand{\bxi}{\bar{\xi}}

\newcommand{\bz}{\bar{z}}

\begin{document}

\bibliographystyle{plain}

\else

\fi

\section{Conclusion}

\subsection{Summary}

With Proposition \ref{consig}, we proved that the signature of a pairing is invariant with regards to cut-and-glue operations. We then proved the other direction for balanced pairings in Theorem \ref{balpair}. Through induction, we extended this result to pairings on surfaces with any number of boundary components in the proof of Theorem \ref{balpairall}. We then applied this result to topological theory to get Corollary \ref{topapp}.

\subsection{Future Directions of Study}

Our work to prove Theorem \ref{balpair} leads us to believe the following conjecture is true:

\begin{conj}\label{allsig} All pairings that have a fixed signature on a surface with connected boundary are connected.

\end{conj}

Just as we generalized Theorem \ref{balpair} to Theorem \ref{balpairall}, one could consider the generalization of conjecture \ref{allsig}:

\begin{conj} All pairings that have a fixed signature on a surface with any number of boundary components are connected.

\end{conj}


\ifx\allfiles\undefined

\bibliography{Index}

\end{document}

\fi



\bibliography{Index}
\end{document}